# The Power and Limits of Predictive Approaches to Observational-Data-Driven Optimization

Dimitris Bertsimas
Massachusetts Institute of Technology, dbertsim@mit.edu, http://web.mit.edu/dbertsim

Nathan Kallus
Cornell Tech and Cornell University, kallus@cornell.edu, http://www.nathankallus.com

While data-driven decision-making is transforming modern operations, most large-scale data is of an observational nature, such as transactional records. These data pose unique challenges in a variety of operational problems posed as stochastic optimization problems, including pricing and inventory management, where one must evaluate the effect of a decision, such as price or order quantity, on an uncertain cost/reward variable, such as demand, based on historical data where decision and outcome may be confounded. Often, the data lacks the features necessary to enable sound assessment of causal effects and/or the strong assumptions necessary may be dubious. Nonetheless, common practice is to assign a decision an objective value equal to the best prediction of cost/reward given the observation of the decision in the data. While in general settings this identification is spurious, for optimization purposes it is only the objective value of the final decision that matters, rather than the validity of any model used to arrive at it. In this paper, we formalize this statement in the case of observational-data-driven optimization and study both the power and limits of predictive approaches to observational-data-driven optimization with a particular focus on pricing. We provide rigorous bounds on optimality gaps of such approaches even when optimal decisions cannot be identified from data. To study potential limits of predictive approaches in real datasets, we develop a new hypothesis test for causal-effect objective optimality. Applying it to interest-rate-setting data, we empirically demonstrate that predictive approaches can be powerful in practice but with some critical limitations.

## 1. Introduction

Data-driven decision-making is transforming modern operations and is being rapidly adopted in practice (McAfee and Brynjolfsson 2012, Brynjolfsson and McElheran 2016). Such data-driven decisions crucially rely on the availability and reliability of very large-scale datasets. But in many decision-making instances, some data is missing and some decision alternatives are not fully characterized by the data. Such issues arise in practice naturally when a manager seeks to choose a decision that performed well historically but the data only contains the realized outcomes of historical decisions and the counterfactual outcomes of potential alternative decisions are – naturally – missing and need to be modeled in order to evaluate the effect of any new managerial decision. These issues can be seen in various data-driven stochastic optimization problems.

There are several prominent examples including ones in revenue management, inventory management, and marketing. In revenue management, a manager may wish to set prices and/or assortments based on a large-scale historical dataset of prices and/or assortments and the demands that resulted from these prices by using demand modeling. Indeed, there are many examples of studies on pricing decisions based on large-scale historical datasets.[1] Besbes et al. (2010) study validating demand models and consider an example of setting loan interest rates based on a large historical dataset of loan offers, including customer acceptance and rejection of each these offers. Ferreira et al. (2016) consider pricing fashion items for an online retailer based on the company's extremely rich historical data on past sale events. And Cohen et al. (2014) study the problem of planning price discounts and consider the application to a large dataset of historical prices for goods in a supermarket and the historical sales generated. Price promotions may also be combined with marketing promotions and similar data-driven decision-making situations arise in these marketing settings. In inventory management, there are examples where demand may be inventory-dependent (Lee et al. 2012).

In all of these examples, the effect of a decision (e.g., price, order quantity, assortment, marketing campaign) on the outcome of an uncertain variable (e.g., demand, yield, returns) has to be modeled from historical, observational data. In general observational data, however, historical decisions and outcomes may be confounded (i.e., have spurious associations), obscuring the isolated, causal effect of any one potential decision. For example, in econometric studies of marketplace supply-demand-price relationships, the analysis generally has to take endogeneity into consideration (Berry et al.

1995, Bijmolt et al. 2005, Phillips et al. 2012). However, in data-driven decision-making, large-scale data is often purely transactional, as in some of the examples above, and may lack the features necessary to enable a precise and sound assessment of the causal effects of any one decision or the strong assumptions necessary for such an assessment may be unfounded. Instead of giving up completely, predictive but not causal analyses are done, modeling the outcome of a given decision as the best prediction of outcome given an observation of the decision in the data. While in general observational settings, where data does not arise from experimental manipulations, this identification is spurious, it may still lead to good decisions.

Indeed, all models are wrong but some models are useful, and, in optimization settings, it is only the objective value of the final decision that matters, rather than the validity of any model used to arrive at it. Thus, while a predictive approach may fail to estimate the causal effect on the objective of the feasible decisions or may be on too shaky ground to provide reliable and valid causal estimates, it may still lead to good decisions down the line. In this paper, we seek to study both the power and limits of predictive approaches to observational-data-driven decision making by studying performance bounds and developing a hypothesis test for causal-effect objective optimality.

To formalize the setting, we consider a decision-making problem given by a stochastic optimization problem. The manager seeks a decision $z \in \mathcal{Z}$. The benefit derived (or negative cost incurred) from the decision $z$ is given by the reward function $r(z)\rho(z,\xi)$, which depends on the realization $\xi$ of an uncertain variable and can be decomposed into a certain part $r(z)$ and an uncertain part $\rho(z,\xi)$. For example, in pricing, $z$ may be price, $\xi$ demand, $r(z) = z - c$ for a procurement cost $c$, and $\rho(z,\xi) = \xi$ so that $r(z)\rho(z,\xi)$ represents net profits. In inventory management, $z$ may be order quantity, $\xi$ demand, $\rho(z,\xi) = -\max\{h(z-\xi), b(z-\xi)\}$ for overage/underage costs, and $r(z) = 1$. The uncertain variable is directly affected by the decision, as demand would be affected by price (or even by order quantity in unique settings as in Lee et al. 2012). Therefore, we consider a collection of random variables, $\xi(z)$ for each $z \in \mathcal{Z}$, representing the random potential outcomes under any one decision. The stochastic optimization problem of interest is

$$z^* = \underset{z \in \mathcal{Z}}{\arg\max}\, \mathbb{E}\left[r(z)\rho(z,\xi(z))\right], \tag{1}$$

**Table 1 Example Demand Data for the Classic MIT Hoodie**

| Day ($i$) | Price ($Z_i$) | Observed Demand ($Y_i$) | Unobserved Demand ($Y_i'$) |
|---|---|---|---|
| 1 | 20 ($) | 1 (units) | 0 |
| 2 | 28 | 0 | 1 |
| 3 | 28 | 1 | 1 |
| 4 | 20 | 2 | 0 |
| $\vdots$ | $\vdots$ | $\vdots$ | $\vdots$ |
| $n-1$ | 20 | 1 | 0 |
| $n$ | 28 | 0 | 1 |

Defining $Y(z) = \rho(z;\xi(z))$ and $y(z) = \mathbb{E}[Y(z)]$, Problem (1) is equivalent to

$$\max_{z\in\mathcal{Z}}\{R(z) := r(z)y(z)\}.$$

In most cases, as in many of the examples reviewed above, the decision is *customized*, e.g., to a specific customer, specific product, or specific instance. For example, Besbes et al. (2010) consider separate interest rates customized for different customer classes defined by ranges of FICO scores, known generally as customized pricing (Phillips 2005). Therefore, we may often think of the distribution in Problem (1) over which expectation is taken as specific to, or conditioned on, a particular customization target, such as customer class. In observational-data-driven decision-making, we consider solving problem (1) based on observational data of $Z$ and $Y = Y(Z)$ (which represents all the relevant information from observing $\xi(Z)$).

For the purpose of illustration, consider a simple example. Table 1 displays the unit demand $Y_i$ observed at the MIT Coop on each day $i = 1,\dots,n$ for the classic MIT hoodie and the price $Z_i$ at which it was offered. Only two prices, \$20 and \$28, have been observed and each was observed $n_{20}$ and $n_{28}$ times, respectively. For each day, there are the demands $Y_i(20)$ and $Y_i(28)$ that would have been observed if the price were set to \$20 or \$28, respectively – these values represent the unseen *demand curve* associated with that day. We only observe $Y_i = Y_i(Z_i)$. We do not observe $Y_i' = Y_i(48 - Z_i)$. For example, on day 1, $Y_1 = Y_1(20) = 1$ and $Y_1' = Y_1(28) = 0$. Using the observed data only, we can compute

$$\tilde{y}_n(20) = \frac{1}{n_{20}}(Y_1 + Y_4 + \cdots + Y_{n-1}) = \text{Average}\left(\{Y_i \,:\, i = 1,\dots,n,\, Z_i = 20\}\right),$$

$$\tilde{y}_n(28) = \frac{1}{n_{28}}(Y_2 + Y_3 + \cdots + Y_n \quad) = \text{Average}\left(\{Y_i \,:\, i = 1,\dots,n,\, Z_i = 28\}\right).$$

If the rows of Table 1 constitute independent and identically distributed (iid) data and $Z$, $Y$, $Y(z)$ represent a generic random draw, then $\tilde{y}_n(20)$ and $\tilde{y}_n(28)$ are our best guesses for the value of demand $Y$ in a new random instance where $Z=20$ or $Z=28$, respectively. In particular, $\tilde{y}_n(z) \to \tilde{y}(z) := \mathbb{E}[Y \mid Z=z]$ almost surely as $n \to \infty$, and $\tilde{y}(z)$ is by definition always the *best* predictor of demand given price (in squared error). Identifying $\tilde{y}(z)$ with the expected demand resulting from a pricing decision $z$ leads to the optimization problem

$$\tilde{z} \in \arg\max_{z \in \mathcal{Z}} \{\tilde{R}(z) := r(z)\tilde{y}(z)\}, \tag{2}$$

where in this particular example $\mathcal{Z} = \{20, 28\}$. When we substitute any estimate $\tilde{y}_n(z)$ of $\tilde{y}(z)$ (or any estimate $\tilde{R}_n(z)$ of $\tilde{R}(z)$) into (2), we call the result a *predictive approach* to pricing from data because it hinges on fitting a *predictive model* of response (or reward) given decision to the data. Given observational data, we can estimate $\tilde{y}(z)$ (or $\tilde{R}_n(z)$) using any regression method, fitting a linear, other parametric, or non-parametric regression to the response $Y$ (or $R(Z)$) to regressor $Z$. For example, Besbes et al. (2010) use non-parametric Nadaraya-Watson kernel regression to construct an estimator $\tilde{R}_n(z)$ of $\tilde{R}(z)$ that is a universally consistent under mild conditions (Greblicki et al. 1984), that is $\tilde{R}_n(z) \to \tilde{R}(z)$ as $n \to \infty$ without specifying a model for $\tilde{R}(z)$.

On average over a new random instance, $y(z) = \mathbb{E}[Y(z)]$ is the expected demand if price were set to $z$ and the problem of selecting a price so to optimize expected profits for a new random instance is given by Problem (1), $\max_{z \in \mathcal{Z}} r(z)y(z)$. In general, $y(z) \neq \tilde{y}(z)$. In particular, if the population distribution of the data is exactly the discrete distribution with weight 1/6 on each of the six displayed rows of Table 1 ($n=6$), then

$$\tilde{y}(20) = 4/3 \neq 7/6 = y(20), \qquad \tilde{y}(28) = 1/3 \neq 1/6 = y(28).$$

This highlights that, in general, Problem (2) is different from Problem (1) for observational data.

Moreover, an estimate of $y(z)$ involves *unobserved* data:

$$y_n(z) = \frac{1}{n}\left(Y_1(z) + Y_2(z) + \cdots + Y_n(z)\right) = \text{Average}\left(\{Y_i(z) \,:\, i = 1, \ldots, n\}\right),$$

In particular, if the data is distributed in the population like the six displayed rows, we can imagine filling in the unseen values in Table 1 with anything, changing $y(z)$, and correspondingly Problem (1), but keeping the observed data completely unchanged. This highlights that, in general, Problem (1) is not necessarily well-specified by the observed data in the most general settings.

However, in an optimization setting, this issue may potentially be moot. Regardless of how close $\tilde{y}(z)$ and $y(z)$ are, it may very well be the case that $\tilde{z}$ is close to $z^*$ and, much more importantly, that $R(\tilde{z})$ is close to $R(z^*)$. For example, in the above example, if the data is distributed in the population like the six displayed rows and supposing the MIT Coop has procurement cost of \$19, we have $\tilde{z} = z^* = 28$ even though $y(z) \neq \tilde{y}(z)$ and $y(z)$ is not well-specified by the observed data, showing the potential power of predictive approaches in practice. But we can also imagine filling in the unseen values in Table 1 such that $\tilde{z} \neq z^*$ and, more importantly, that $R(z^*)$ is much larger than $R(\tilde{z})$, showing the potential limits of predictive approaches in practice.

In this paper, we provide a thorough exploration of this issue in data-driven decision-making, prove new performance guarantees in certain common decision-making situations, and develop a new hypothesis test to evaluate objective optimality based on observational data. We explore the gap between Problems (1) and (2) and review the issue of identifiability and confounding (Section 2) and show through examples the implications in the context of data-driven decision-making. In the specific instance of price optimization problems, we bound the suboptimality of any predictive approach with respect to the true optimum even when the optimum cannot be identified from the data (Section 3). Our bounds leverage the special structure of the optimization problem as well as special features that are common in the data usually available in practice. This demonstrates the power of predictive approaches to observational-data-driven optimization. To study potential limits of predictive approaches in real datasets, we develop a hypothesis test for the objective optimality of any data-driven solution in Problem (1), such as solutions based on observational data (Section 4). The test allows us to determine whether the reward generated by any one decision-making algorithm, such as a predictive one, can be distinguished as suboptimal to a statistically significant degree based on purely observational data. To develop the test, we show favorable asymptotics (consistency and asymptotic normality) for a non-parametric solution that works under certain identifiability conditions (Section 4.2) and use this to establish an asymptotic null distribution for a test statistic for the hypothesis of objective optimality (Section 4.3). We also present a parametric approach to solving Problem (1) under these conditions where the solution is identifiable (Section 5). Using our hypothesis test, we empirically study an interest-rate-setting problem with data from Columbia University Center for Pricing and Revenue Management (2012) (Section 6). We demonstrate mixed results showing that predictive approaches are sometimes distinguishable as suboptimal to a statistically significant and sometimes not, depending on the setting. This finding shows that the distinction between Problems (1) and (2) is of real, practical relevance and that predictive approaches can be powerful in practice but with some limits. We are also able to demonstrate that a parametric approach that correctly accounts for the observational data (under the appropriate identifiability conditions) can recover 36-70% of profits lost by predictive approaches. This second finding expands the scope of recent work on the sufficiency of parametric models for revenue management (Besbes et al. 2010, Besbes and Zeevi 2015) to revenue management from observational data.

## 2. Identifiability and the Gap Between Problems (1) and (2)

As exemplified in the introduction, the mean-response curve $y(z)$ that defines Problem (1) is not well-specified by observational data and is distinct from $\tilde{y}(z)$, which is. In this section, we first consider a more detailed example and study its implications on identifiability of $z^*$ (which we will define precisely) and then we consider the gap between $y(z)$ and $\tilde{y}(z)$ in preparation for bounding the gap between $R(z)$ and $R(\tilde{z})$ in the next section.

Example 1 (Consulting for the MIT Coop). Alice and Bob are hired by the MIT Coop to help determine an optimal sale price for the classic MIT hoodie, which the MIT Coop procures at a unit price of $c = \$19$ ($r(z) = z - 19$). The MIT Coop is debating between a retail price of \$20 and a retail price of \$28. In any given day in the past, the MIT Coop has offered the hoodie at

**Table 2 Data for Example 1**

| $Z$ | $Y$ | $\mathbb{P}$ |
|---|---|---|
| 20 | 0 | 0 |
| 20 | 1 | 8/18 |
| 20 | 2 | 1/18 |
| 28 | 0 | 8/28 |
| 28 | 1 | 1/18 |
| 28 | 2 | 0 |

(a) Joint distribution of historical price and demand

| $Z$ | $Y(20)$ | $Y(28)$ | $\mathbb{P}$ |
|---|---|---|---|
| 20 | 1 | 0 | 32/81 |
| 20 | 1 | 1 | 4/81 |
| 20 | 2 | 0 | 4/81 |
| 20 | 2 | 1 | 1/162 |
| 28 | 1 | 0 | 32/81 |
| 28 | 1 | 1 | 4/81 |
| 28 | 2 | 0 | 4/81 |
| 28 | 2 | 1 | 1/162 |

(b) Alice's demand model

| $Z$ | $Y(20)$ | $Y(28)$ | $\mathbb{P}$ |
|---|---|---|---|
| 20 | 1 | 0 | 40/99 |
| 20 | 1 | 1 | 4/99 |
| 20 | 2 | 0 | 0 |
| 20 | 2 | 1 | 1/18 |
| 28 | 1 | 0 | 4/9 |
| 28 | 1 | 1 | 2/45 |
| 28 | 2 | 0 | 0 |
| 28 | 2 | 1 | 1/90 |

(c) Bob's demand model

either of the two prices and observed either no units sold, one unit sold, or two units sold. The Coop has a great deal of observational data.

Alice and Bob collate the data into a table that shows the frequency of each price-demand combination over history shown in Table 2(a). Due to the abundance of data, Alice and Bob are confident that this is a faithful representation of the joint distribution of $(Y, Z)$. Naturally, the data only has the demand that was in fact observed and the demand that would have been observed under any other price is missing. A full demand model models the distribution of the demand curve $Y(\cdot)$ for a new sale event.

Alice regresses demand on price by computing a weighted average in each of the columns of Table 2(a) and finds that $\mathbb{E}\left[Y \middle| Z=20\right] = 10/9$, $\mathbb{E}\left[Y \middle| Z=28\right] = 1/9$. She constructs a demand model wherein $y(z) = \mathbb{E}\left[Y \middle| Z=z\right]$ and arrives at the one shown in Table 2(b). Alice verifies that her model fully agrees with the observed data (via the transformation $Y = Y(Z)$) and computes $R(20) = (20-19) \times (\frac{72}{81} \times 1 + \frac{9}{81} \times 2) = 10/9$ and $R(28) = (28-19) \times (\frac{72}{81} \times 0 + \frac{9}{81} \times 1) = 1$, concluding that $z^\star = 20$ is the optimal price.

Bob, working from home that day and unaware of Alice's progress, has independently come up with another model, shown in Table 2(c), in order to explain the observed pricing data. Bob, too, verifies that his model completely agrees with the observed data and calculates $R(20) = (20-19) \times (\frac{14}{15} \times 1 + \frac{1}{15} \times 2) = 16/15$ and $R(28) = (28-19) \times (\frac{28}{33} \times 0 + \frac{5}{33} \times 1) = 15/11$, concluding, differently from Alice, that $z^\star = 28$ is in fact the optimal price.

Alice and Bob had both come up with demand models that fully concur with the observed data but recommended different prices as optimal. Both models support the data fully. Both models, as well as the data, fully agree with a homoscedastic linear model,

$$Y = \frac{65}{18} - \frac{1}{8}Z + \omega, \qquad \omega = \begin{cases} -1/9, & \text{with prob. } 8/9, \\ \ \ 8/9, & \text{with prob. } 1/9, \end{cases} \qquad \omega \perp\!\!\!\perp Z,$$

where regressor $Z$ is independent of zero-mean error $\omega$. Since the two models agree on this form but recommended different prices, this highlights that this is not an issue of misspecifying a functional form for the demand model.

### 2.1. Identifiability

The issue we encountered above is one of identifiability and shows that $z^*$ is non-identifiable in general.

Definition 1. Let $\Pi = \{\mathbb{P}_\theta : \theta \in \Theta\}$ be a model for the distribution of the observed data. We say that $\phi : \Theta \to \Phi$ is *identifiable* if for any $\mathbb{P}_{\theta_1}, \mathbb{P}_{\theta_2} \in \Pi$ such that $\mathbb{P}_{\theta_1} = \mathbb{P}_{\theta_2}$, we have $\phi(\theta_1) = \phi(\theta_2)$.

In the above definition, $\Theta$ and $\Phi$ may be arbitrary sets, that is, the model need not be parametric. Note that if any $\phi$ is not identifiable then any finer quantity, such as $\theta$ itself, is not identifiable.

To connect the above definition with our decision-making setting, we let $\theta$ denote the joint distribution of $(Z, Y(\cdot))$, we let $\mathbb{P}_\theta$ be the corresponding distribution of the data $(Z, Y)$ (begotten via the transformation $Y = Y(Z)$), and we let $\phi$ map $\theta$ to the optimal decision $z^\star$ (or, set thereof) as described by Problem (1). Then, Example 1 above proves the following result:

Corollary 1. *The optimal decision $z^\star$ is not identifiable on the basis of observations of $(Y, Z)$.*

In fact, we have shown a stronger result:

Theorem 1. *The optimal decision $z^\star$ is not identifiable on the basis of observations of $(Y, Z)$ even under the Gauss-Markov assumptions:*

*a. Linearity: there is a random variable $\omega$ such that $Y = \beta_0 + \beta_1 Z + \omega$.*

*b. Exogeneity of independent variables: $\mathbb{E}\left[\omega \middle| Z\right] = 0$.*

*c. Homoscedasticity: $\mathrm{Var}\left(\omega \middle| Z\right) = \mathrm{Var}\left(\omega\right)$ is constant.*

*d. No collinearity: $Z$ is not constant.*

In Example 1, exogeneity and homoskendasticity are a consequence of $\mathbb{E}[\omega] = 0$ and $\omega \perp\!\!\!\perp Z$. Exogeneity implies $\mathrm{Cov}(\omega, Z) = 0$. Note that whenever the optimal price is not identifiable, the mean-response $y(z)$, a finer quantity, is not identifiable either.

When *is* $z^*$ identifiable from observational data? Naturally, one case is when $y(z) = \tilde{y}(z)$, since $\tilde{y}(z)$ is always identifiable from observational data. Next, we study the gap between $y(z)$ and $\tilde{y}(z)$.

### 2.2. The Confounding Error Gap

Let $\epsilon(z) = Y(z) - y(z)$ be the deviation of the response curve from its mean. Following our convention for $Y = Y(Z)$, we also let $\epsilon = \epsilon(Z) = Y - y(Z)$, which is the deviation of the *observed* $Y$ under decision $Z$ from the mean response that would be induced by the decision $Z$. The degree to which $\epsilon$ is correlated with $Z$ is known as *confounding*. In fact it can be directly related to the discrepancy between $y(z)$ and $\tilde{y}(z)$,

$$E(z) := \tilde{y}(z) - y(z) = \mathbb{E}[Y \mid Z = z] - y(z) = \mathbb{E}[Y - y(Z) \mid Z = z] = \mathbb{E}[\epsilon \mid Z = z].$$

We call $E = E(Z) = \mathbb{E}[\epsilon \mid Z]$ the confounding error.

Thus, the error $E$ is directly related to the association of historical decisions $Z$ and the unique idiosyncrasies $\epsilon$ of historical outcome events. Independence of the two – i.e., that $Z$ is independent of the particular event, as in a randomized controlled trial (RCT) – would immediately imply that the confounding error is exactly zero: $E = \mathbb{E}[\epsilon \mid Z] = \mathbb{E}[\epsilon] = \mathbb{E}[Y(Z)] - \mathbb{E}[y(Z)] = 0$. Correspondingly, we would have $\tilde{y}(z) = y(z)$. Therefore, when the data is the result of experimental manipulations rather than observation, we must have $\tilde{z} = z^*$. In the customized setting, all distributions are conditional on the customization target, and the independence necessary for this is the conditional independence given this target. We use a related condition in Section 4.1.2.

Note the critical distinction between $\epsilon$ and the regression errors (residuals) $\omega = Y - \mathbb{E}\left[Y \middle| Z = z\right]$. Regression errors, by their very definition, will *always* have $\mathbb{E}\left[\omega \middle| Z\right] = \mathbb{E}[Y \mid Z] - \mathbb{E}[Y \mid Z] = 0$ even in the most general setting for observational data, but the same is generally only true of $\epsilon$ in experimental settings. Thus, in general settings, $\mathbb{E}[\epsilon \mid Z] \neq 0$ and $\tilde{y}(z) \neq y(z)$.

**Figure 1 Predictive demand and revenue under confounding.**

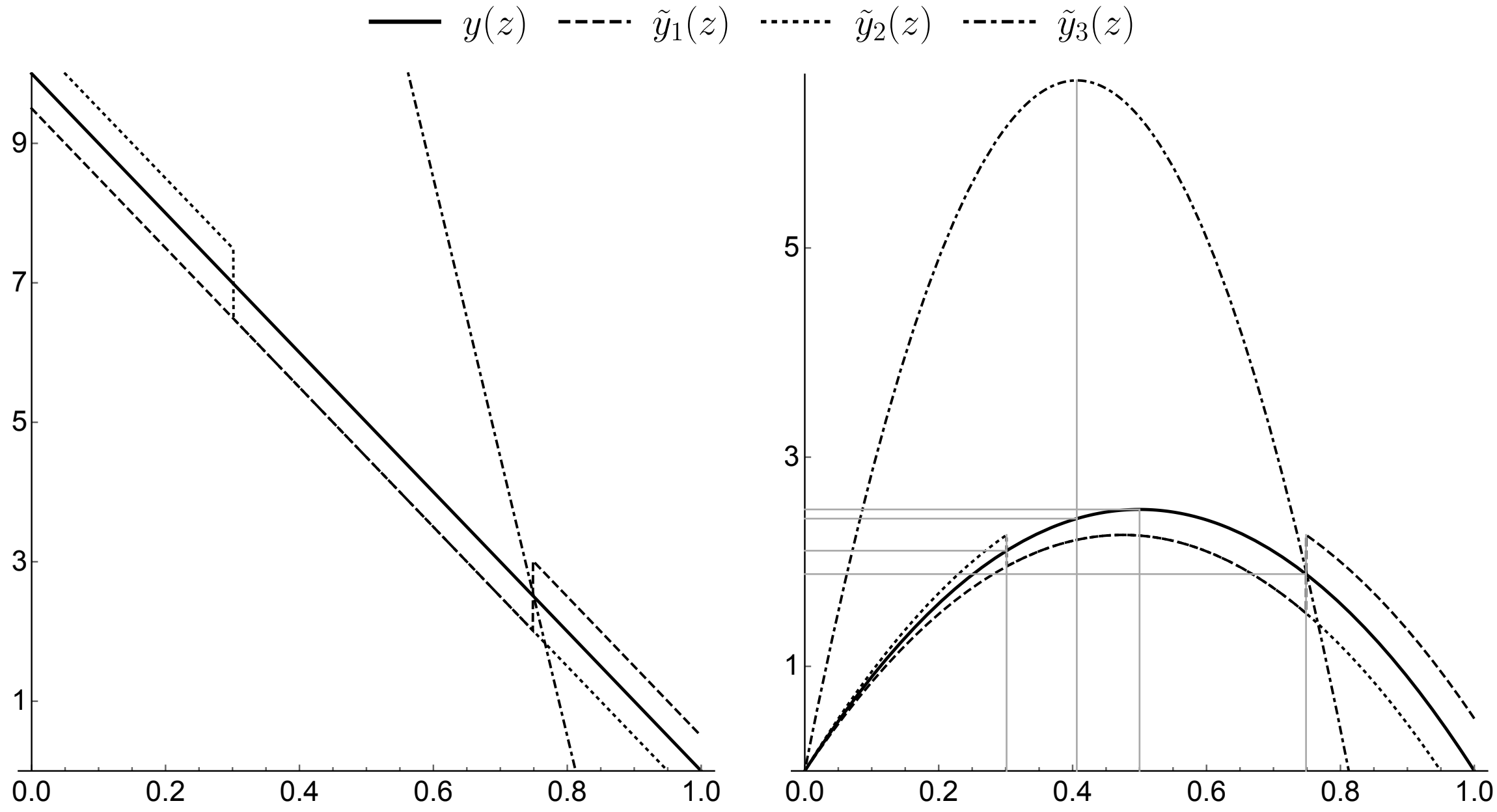

(a) Demand

(b) Revenue

However, even when $\tilde{y}(z)$ is very different from $y(z)$, what matters really is the performance of the resulting solution, $R(\tilde{z})$. Figure 1 shows an example for the pricing problem with $r(z) = z$ and $\mathcal{Z} = [0,1]$. Figure 1(a) shows one true curve for $y(z)$ and three possible different curves for $\tilde{y}(z)$. The first two, $\tilde{y}_1(z)$ and $\tilde{y}_2(z)$, differ from $y(z)$ by a constant amount 1/2 throughout the domain and the second is non-increasing. The third, $\tilde{y}_3(z)$, differs from $y(z)$ a lot in magnitude but not in its linear shape. Figure 1(b) shows the corresponding revenue curves. Optimizing each $r(z)\tilde{y}_i(z)$ to find $\tilde{z}_i$ and plugging each into $R(z)$, we see that $R(\tilde{z}_1)/R(z^\star) = 0.75$, $R(\tilde{z}_2)/R(z^\star) = 0.84$, $R(\tilde{z}_3)/R(z^\star) = 0.96$. This shows a few examples of how discrepancies between $\tilde{y}(z)$ and $y(z)$ translate to performance. Since in general settings $z^\star$ can never be known, to have performance that is reasonably within that of the unknown and optimal $z^\star$ would be very fortunate. An important question we address next is how to bound this optimality gap in general, without knowledge of $y(z)$, so that we can theoretically guarantee good performance even if the best decision cannot be pinned down.

## 3. The Power of Predictive Approaches to Pricing

In this section, we focus in on the pricing instance of Problem (1). In particular, we have that $r(z) = z - c$ for some procurement cost $c > 0$, $\mathcal{Z} = (c, \infty)$ is univariate, $Y(z)$ represents the *random* demand curve, and $y(z)$ represents the average demand curve.

To show the power of predictive approaches, we establish bounds on how suboptimal a predictive approach may be in the pricing problem. Such bounds show that even if confounding is present and even if Problem (1) is not well-defined given the data, predictive approaches, which are implementable in practice, can still get us reasonably close to the optimum under certain assumptions. The bounds leverage both the special structure of the pricing problem and of pricing data in practice and address error in the metric of interest, which is true profit $R(p)$. The intent is to express this error in *relative* terms, using only *few* parameters, and using assumptions we can *reason* about. Beyond pricing, these bounds suggest a framework for establishing similar results in other observational-data-driven decision-making problems by appealing to the special structure of the problem and not just the error in misspecification.

The bounds we present are expressed in terms relative to the *size of the market*:

$$y_0 = \sup_{z \in \mathcal{Z}} y(z).$$

The bound is based on the *magnitude* of confounding error relative to the size of the market.

THEOREM 2. *In the pricing problem, if $y(z)$ is non-increasing and linear and $|E|/y_0 \le \gamma$ then*

$$1 - 4\gamma - 4\gamma^{3/2} - \gamma^2 \le \frac{R(\tilde{z})}{R(z^\star)} \le 1.$$

The proof is given in the appendix. The bound is non-negative for values of $\gamma$ up to $3 - \sqrt{8} \approx 0.17$. We plot the bound in Figure 2(a). Applying this bound to $\tilde{y}_1(z)$ of Figure 1, which has $\gamma = 0.05$, we get 0.71, which indeed bounds the observed ratio of 0.75.

In general, the magnitude of the confounding error is neither known nor estimable from data. Our next bound seeks to leverage particular structure that is symptomatic of pricing data to express suboptimality in terms of average demand, which can be estimated from data.

**Figure 2 The Bounds of Theorems 2, 3, and 4.**

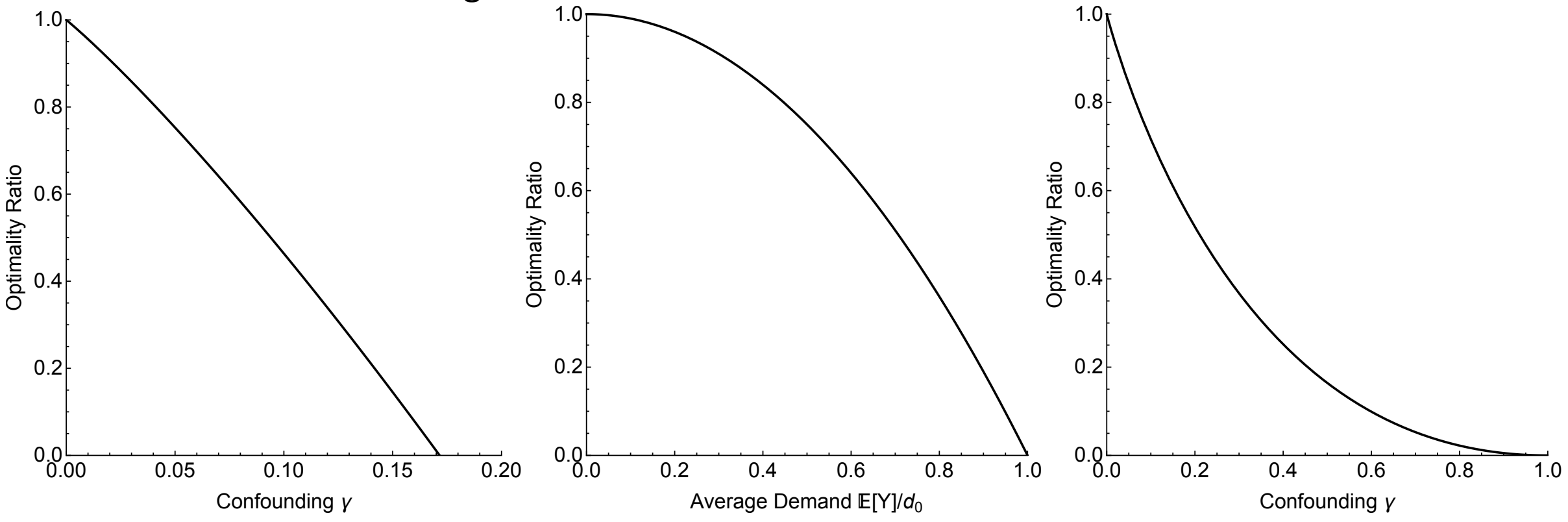

(a) Theorem 2 (b) Theorem 3 (c) Theorem 4

Theorem 3. *In the pricing problem, if $y(z)$ is linear and non-increasing and $\epsilon$ and $Z$ are non-positively correlated and jointly normal then*

$$1-\left(\frac{\mathbb{E}\left[Y\right]}{y_0}\right)^2 \le \frac{R(\tilde{z})}{R(z^\star)} \le 1.$$

The proof is given in the appendix. Note that $\mathbb{E}\left[D\right]$ can be unbiasedly and consistently estimated by the sample average of observed demands. The assumption of non-positive correlation between $\epsilon$ and $P$ corresponds to a common feature of pricing datasets. For example, promoting a product via advertising, which would increase its potential demand at any given price, would often coincide with promotion via price discounts and this would lead to such non-positive correlation. This assumption, which can be reasoned about in such a way, leads to a stronger bound that is completely independent of the size of confounding error. We plot the bound in Figure 2(b).

The curve $\tilde{y}_3(z)$ of Figure 1 can be seen to arise from a situation where $y(z)$ is as in Figure 1, $Z$ and $\epsilon$ are jointly normal with a covariance that is $-30$ of the variance of $Z$, and $\mathbb{E}\left[Y\right] = 2.5$. Applying Theorem 3, we get 0.94, which indeed bounds the observed ratio of 0.96. (In particular, we can achieve the bound by letting the covariance approach $-\infty$.)

If we consider a similar non-positive relationship between $\epsilon$ and $Z$ in the setting of Theorem 2, we obtain an improved bound.

THEOREM 4. *In the pricing problem, if $y(z)$ is non-increasing and linear, $E(z)$ is non-increasing, and $|E|/y_0 \le \gamma \le 1$ then*

$$1-4\gamma+4\gamma^{3/2}-\gamma^2 \le \frac{R(\tilde{z})}{R(z^\star)} \le 1.$$

The proof is given in the appendix. This bound is similar to that in Theorem 2 but the additional assumption that $E(z)$ is decreasing allows us to achieve a strictly stronger bound that is non-negative for values of $\gamma$ up to 1. That $E(z)$ is non-increasing captures the same common feature of pricing datasets in a more general way. We plot the bound in Figure 2(c). Applying this bound to $\tilde{y}_2(z)$ of Figure 1, which has $\gamma = 0.05$ and $E(z)$ decreasing, we get 0.79, which indeed bounds the observed ratio of 0.84.

## 4. A Test for Causal-Effect Optimality

In the last sections we saw that predictive approaches, while in fact solving the wrong problem, can often lead to good objective performance nonetheless. But, naturally, there are limits to their power. To evaluate the performance of predictive approaches more generally and check whether they can actually be distinguished from optimal in practice, we next develop a hypothesis test for causal-effect optimality. We do this in the spirit of Besbes et al. (2010), who develop a test for objective optimality in *experimental* settings to test whether parametric models suffice to achieve good objective performance. Our test extends the work to the case of observational data and to causal effects.

Each in their own context, Besbes et al. (2010), Cohen et al. (2014), Ferreira et al. (2016) address an observational-data-driven revenue management problem by first estimating a predictive model of (usually, customized) demand and then optimizing a price or other management decision. It may be argued that it is important that estimated profits faithfully represent true profits in the resulting optimization problem, but in fact this point is moot insofar as actual profits generated by the resulting pricing strategy are satisfactory. Moreover, if a fully faithful model is not feasible for lack of identifiability, a predictive approach is all that may be hoped for.

Suppose we wish to test the objective optimality in Problem (1) of a data-driven decision-making algorithm that prescribes the decision $\hat{z}_n$ based on $n$ data points. Let $\hat{z}$ be the corresponding *full-information* decision algorithm, i.e., the hypothetical decision this algorithm would pick given infinite data $n \to \infty$. For complete generality, we leave the meaning of this vague and only require the following condition in defining what $\hat{p}$ means.

Assumption 1 **(Convergent Decision-Making).** $\hat{z}_n - \hat{z} = O_p(1/\sqrt{n})$ *for some fixed* $\hat{z} \in \mathcal{Z}$.

The notation $U_n = O_p(a_n)$ means that for any $\epsilon > 0$ there is $M > 0$ such that $\mathbb{P}(\|U_n/a_n\| > M) < \epsilon$ eventually. In particular, if $U_n/a_n$ converges in distribution (to anything) then $U_n = O_p(a_n)$.

For example, the strategy $\hat{z}_n$ that fits a Nadaraya-Watson kernel regression to estimate conditional $\tilde{R}(z)$ and selects $\hat{z}_n$ by plugging the estimate into Problem (2) will (under some regularity) eventually arrive at the decision $\hat{z} = \tilde{z}$ that optimizes the hypothetical Problem (2), and, in particular, satisfies $\hat{z}_n - \hat{z} = O_p(1/\sqrt{n})$ because the difference is asymptotically normal (Ziegler 2002). A similar condition, with a potentially different $\hat{z}$, is true of the strategy that uses a parametric maximum-likelihood regression (Besbes et al. 2010, cf. Lemma 2 and Lipschitz condition in the proof of Lemma 3). Assumption 1 is generally true of most data-driven approaches.

We would like to test whether the nominal decision that our data-driven decision-making algorithm is really getting at, that is, $\hat{z}$, is truly optimal or not. Therefore, we would like to test the following null hypothesis $H_0$ against the alternative $H_1$:

$$H_0 : R(z^\star) = R(\hat{z}), \quad H_1 : R(z^\star) > R(\hat{z}).$$

That is, we would like to test whether the nominal price that our pricing strategy would be prescribing is generating optimal profits.[2]

A test for the hypothesis $H_0$ can be interpreted as rejecting a pricing strategy if it *generates profits that are distinguishable from optimal to a statistically significant degree based on the data.*

### 4.1. The Setting for Testing $H_0$

We can repeat the argument in Section 2.1 to see that in the most general settings the truth value of $H_0$ is not identifiable from observational data. Therefore, we can have little hope of testing it in the most general settings.

Instead, we consider a situation where additional variables are available to account for all confounding factors and ensure the identifiability of $H_0$. What this allows us to do is evaluate the hypothesis $H_0$ in real, practical settings where such data is available so to understand whether predictive approaches, that ignore this data, work, with the hope of implying whether they work in general real, practical settings where such additional data may or may not be available.

**4.1.1. The On-Line Auto Lending Case.** An exceedingly appropriate case study is given by the on-line auto lending data of Columbia University Center for Pricing and Revenue Management (2012) and the customized rate-setting problem considered by Besbes et al. (2010). The on-line auto lending data consists of past sale events where a customer fills out a loan application, if approved an interest rate is quoted (where rate is essentially equivalent to a price), and the customer either accepts or rejects the loan (binary demand). Besbes et al. (2010) study the problem of prescribing interest rates for automobile loans based on this data and customized to each of eight customer segments delineated by predefined ranges of FICO scores, term lengths, and season (see Example 4 for additional detail). The authors estimate predictive models for $\tilde{y}(z)$ or $\tilde{R}(z)$ within each segment by separately regressing $Y$ or $R(Z)$ on $Z$ based on the data available from each segment. The focus of their study is evaluating whether parametric regression (specifically logistic regression) methods lead to interest rates that are optimal in the predictive Problem (2) by evaluating the rates in a non-parametric (Nadaraya-Watson kernel) regression estimate of $\tilde{R}(z)$.

The dataset description says that approval and rate is based on "credit information and other criteria." Such criteria would almost certainly also be associated with the potential likelihood of the consumer to accept a loan offer at any one particular rate (the demand curve). Even if rates are not chosen strategically in response to demand, they could be chosen based on default risk or expected loss, which may be in turn associated with the demand curve. Therefore, we argue that confounding is likely present, i.e., $\mathbb{E}[\epsilon \mid Z] \neq 0$ and $y(z) \neq \tilde{y}(z)$. Moreover, even within each of the eight customer segments considered by Besbes et al. (2010), we would argue that confounding is likely present conditioned on segment because the FICO ranges considered are predefined (rather than recursively segmented to fully condition on FICO score) and FICO score only account for credit information and not all of the "other criteria" that influence historical rates. Therefore, since the non-identifiability caused by confounding is not a model-specification issue, even the non-parametric estimates may estimate only the predictive objective $\tilde{R}(z)$. The question is whether the resulting interest rates perform well in the true objective $R(z)$.

The on-line auto lending data contains much more information about each loan applicant and the associated sale event, including the precise FICO credit score of the applicant, the length of the term over which the loan is to be repaid, the dollar amount of the loan, whether the car to be purchased is new, used, or refinanced, competitors' rate, prime rate, and who referred the applicant. From here on, we let $X$ denote these covariates. If the covariates $X$ encompass the aforementioned "credit information and other criteria," then we may in fact be able to identify $y(z)$ and test $H_0$, checking the objective optimality of the predictive approaches above, whether parametric or non-parametric. Note that the customization target, a variable taking values in $\{1,\dots,8\}$, is much coarser than $X$, i.e., it is a highly non-injective function of $X$. As a finite-valued function of $X$, it suffices for us to consider the data, including $X$, within each segment separately.

**4.1.2. Identifiability of $H_0$.** Let $X$ denote some covariates observed concurrently with each historical outcome event, as in the example of the autoloan dataset, and let $X_i$ be the observation corresponding to the $i^{\text{th}}$ event. For example, in a pricing problem, the covariates $X$ may include, for example, characteristics of the customer, whether a product was featured in a promotional flyer, external signals about demand used for pricing, etc. Even with this data, for the moment, we still restrict ourselves to the problem of choosing a single price for the whole population of sale events, or customized to a discrete segment.

Sometimes covariates $X$ can help us disassociate the random variable $Z$ and the particular response curve $Y(z)$, the association between which is the source of confounding. One such sufficient condition for $X$ to account for all such association is the following standard ignorability condition (Hirano and Imbens 2004), a continuous version of ignorability in Rosenbaum and Rubin (1983).

Assumption 2. *For every $z \in \mathcal{Z}$, we have that, conditioned on $X$, $Y(z)$ is independent of $Z$. That is, $\forall z \in \mathcal{Z},\ Y(z) \perp\!\!\!\perp Z \,\big|\, X$.*

The condition says that, historically, $X$ accounts for all the event-specific features that may have influenced the decision $Z$ up to idiosyncratic and independent randomness in $Z$. For example, we argue that $X$ in the on-line auto lending dataset satisfies this condition.

Under this condition, it is immediate and well-understood that we have identifiability:

$$y(z) = \mathbb{E}\left[\mathbb{E}\left[Y(z) \mid X\right]\right] = \mathbb{E}\left[\mathbb{E}\left[Y(z) \mid Z = z, X\right]\right] = \mathbb{E}\left[\mathbb{E}\left[Y(Z) \mid Z = z, X\right]\right] = \mathbb{E}\left[\mathbb{E}\left[Y \mid Z = z, X\right]\right], \quad (3)$$

where the first equality is by iterated expectations, the second is by Assumption 2, the third is by $Z = z$, and the last by $Y = Y(Z)$. The last expectation is expressed solely in terms of the joint distribution of $(Z, X, Y)$, which gives the identifiability $y(z)$ and hence the objective function $R(z)$. The optimal decision is given by optimizing $R(z)$ and the truth value of $H_0$ is given by plugging in values to the objective function $R(z)$. So these are also identifiable. We summarize this as follows.

Theorem 5. *Under Assumption 2, the mean-response curve $y(z)$, optimal decision $z^\star$, and the truth value of $H_0$ are all identifiable on the basis of observations of $(X, Y, Z)$.*

Note the last expectation in equation (3), $\mathbb{E}\left[\mathbb{E}\left[Y \middle| Z = z, X\right]\right]$, is not conditioned on $Z = z$ and cannot be marginalized via iterated expectations. In words, it says to take the average of the conditional expectation of $Y$ given $X = x$, $Z = z$ over all $x$ using the *marginal* distribution of $X$.

If $Z$ were chosen without regard to any specific event, Assumption 2 holds with a null $X$ variable (formally, $\sigma(X) = \{\Omega, \varnothing\}$). In particular, this is the case in dynamic demand learning and pricing as in Bertsimas and Perakis (2006), Besbes and Zeevi (2009), Harrison et al. (2012) because each sale event is assumed independent and nothing about a present sale event is considered when setting the price. This is the experimental setting where $Z$ is the result of controlled, randomized manipulation, which is rarely the case in practice for observational data. If there is not sufficient recorded information in $X$ to merit Assumption 2, it is said that there is residual endogeneity. In this case, Theorem 5 fails, but there may be other conditions that enable identification such as the availability of instrumental variables (see e.g. Bijmolt et al. 2005). Here we focus on the case where Assumption 2 holds.

**Table 3 Data for Example 2**

| | $X=0$ | | $X=1$ | |
|---|---|---|---|---|
| | $P=20$ | $P=28$ | $P=20$ | $P=28$ |
| D=0 | 0 | 4/9 | 0 | 0 |
| D=1 | 4/9 | 2/45 | 0 | 1/90 |
| D=2 | 0 | 0 | 1/18 | 0 |

Let us consider Assumption 2 and its ramifications in an example.

Example 2 (Consulting for the MIT Coop). Consider again the hypothetical case of Example 1. Recall, Alice and Bob both came up with models for demand that completely agreed with the data but gave rise to different optimal prices. Thus, we concluded that the data observed could not possibly identify the right optimal price.

Suppose Assumption 2 holds with $X$ being a null variable, i.e., without any extra information. This condition eliminates Bob's model – it no longer agrees with both the data and this condition. On the other hand, Alice's model remains valid – in fact it turns out to be the unique model that agrees with both the data and this condition. Hence, under this condition, $z^\star = 20$ is the correct optimal price. But for Assumption 2 to hold with $X$ being a null variable we would have needed experimental data, where prices are set at random for the sake of experiment.

Suppose instead that we recorded additional information about each sale event: whether there was a major home game that day ($X = 1$) or not ($X = 0$). On average, there is a game 2 days of each month ($\mathbb{P}(X = 1) = 2/30$). Suppose tallying the historical observations led to the summary of the data shown in Table 3. If prices were chosen independently of whether there was a game, the previous scenario still holds and Alice's model is the uniquely correct one. If prices were more often cut to \$20 when there was a game, then we are no longer in the experimental setting. In fact, if we assume Assumption 2 holds with this $X$, then it turns out that Alice's model is ruled out and Bob's model is the unique model that accommodates both this condition and the data observed, in which case $z^\star = 28$ is the correct optimal price. In this hypothetical example, we are seeking a universal price, to be set a priori without regard to whether there is a game, but the price can also be customized.

With $H_0$ identifiable, what remains is to construct a hypothesis test. To do so, we first construct non-parametric estimates for $R(z)$ and $z^\star$, which will allow us to construct a test statistic.

### 4.2. Non-Parametric Estimates of $R(z)$ and $z^\star$

In this section, we study non-parametric estimates of $R(z)$ and $z^\star$ under Assumption 2 that are model-independent in that they will converge to the true objective and optimal decision regardless of the true underlying distribution, given sufficient data and some regularity assumptions. Henceforth, we assume that $X$ is a vector of covariates taking values in $\mathbb{R}^k$ and that $\mathcal{Z} \subseteq \mathbb{R}^\ell$.

The proof of Theorem 5 says that, under Assumption 2, the objective function can be written as $R(z) = \mathbb{E}\left[\mathbb{E}\left[r(Z)Y \middle| Z = z, X\right]\right]$. Thus, to estimate it, one approach may be to estimate the regression function $\mathbb{E}\left[r(Z)Y \middle| Z = z, X = x\right]$ and then average the estimated function over an estimate for the marginal distribution of $X$. Then, the optimizer of this estimate can be used as an estimator for the optimal decision.

First, one non-parametric estimate of the marginal distribution of $X$ is simply the empirical distribution, which places unit mass at each of the observations $X_i$. Second, to estimate the regression function $\mathbb{E}\left[r(Z)Y \middle| Z = z, X = x\right]$ non-parametrically, we can use Nadaraya-Watson kernel regression (Nadaraya 1964, Watson 1964). The estimate, based on a choice of kernel $K : \mathbb{R}^{\ell+k} \to \mathbb{R}$ and bandwidth $h_n > 0$, is

$$\overline{R}_n(z, x) = \frac{\sum_{i=1}^n K(\frac{z-Z_i}{h_n}, \frac{x-X_i}{h_n}) r(Z_i) Y_i}{\sum_{i=1}^n K(\frac{z-Z_i}{h_n}, \frac{x-X_i}{h_n})}, \tag{4}$$

where $K(\frac{z-Z_i}{h_n}, \frac{x-X_i}{h_n}) = K(\frac{z_1-Z_{i1}}{h_n}, \dots, \frac{z_k-Z_{ik}}{h_n}, \frac{x_1-X_{i1}}{h_n}, \dots, \frac{x_k-X_{ik}}{h_n})$. This regression estimator arises as the conditional expectation with respect to the Parzen window density estimator (Parzen 1962) for the joint density of $(Z, X, r(Z)Y)$. A kernel function can be any function that has positive, finite integral. A kernel mimics a continuous distribution centered at the data points, the width of which is determined by the bandwidth. There are a variety of kernels used in practice (Härdle 1990).

Combining the two estimators as detailed above, we arrive at the following estimate for the objective function

$$\overline{R}_n(z) = \frac{1}{n}\sum_{i=1}^{n} \overline{R}_n(z, X_i) = \frac{1}{n}\sum_{i=1}^{n} \frac{\sum_{j=1}^{n} K\left(\frac{z-Z_j}{h_n}, \frac{X_i-X_j}{h_n}\right) r(Z_j)Y_j}{\sum_{i=1}^{n} K\left(\frac{z-Z_i}{h_n}, \frac{Z_i-Z_j}{h_n}\right)}. \tag{5}$$

Optimizing the above estimate over $\mathcal{Z}$ yields a non-parametric observational-data-driven decision-making algorithm

$$\overline{z}_n \in \arg\max_{z\in\mathcal{Z}} \overline{R}_n(z). \tag{6}$$

One question that arises is how do these decisions behave asymptotically. In particular, does this algorithm lead to a decision and objective value that converge to the optimal decision and objective value. Since the estimates are non-parametric, the hope is that this can occur under model-free assumptions. Next we show that this is indeed the case under the following regularity conditions, focusing on the univariate case $\mathcal{Z} \subseteq \mathbb{R}$.

Assumption 3 **(Kernel Conditions)**.

*a.* $0 < \int_{\mathbb{R}^{1+k}} K(u)\, du < \infty$.

*b.* *$K$ is zero outside a bounded set.*

*c.* *$K$ is twice Lipschitz-continuously differentiable.*

*d.* *$K$ has order at least $s \in \mathbb{N}$, that is, $\int K(u)u^\alpha du = 0 \quad \forall \alpha \in \mathbb{N}^{1+k} : |\alpha| < s$.*

*e.* *$h_n \to 0$ and $nh_n^{2s+3} \to 0$.*

*f.* *$nh_n^{k+5}/\log(n) \to \infty$ and $nh_n^{2k+1}/\log(n)^2 \to \infty$.*

Assumption 4 **(Optimality Conditions)**.

*a.* *$\mathcal{Z} \subseteq \mathbb{R}$ is compact.*

*b.* *$z^\star$ uniquely maximizes $R(z)$ on $\mathcal{Z}$.*

*c.* *$z^\star$ lies in the interior of $\mathcal{Z}$.*

*d.* *$R(z)$ is twice continuously differentiable and $R''(z^\star) < 0$.*

Assumption 5 **(Distributional Conditions)**.

*a.* *$X$ and $Z$ are continuously distributed on a compact support where the joint density, $f_{Z,X}(z,x)$, is bounded away from zero.*

*b.* *The marginal density of $X$, $f_X(x)$, is bounded and continuously differentiable.*

*c.* $\mathbb{E}[Y^4] < \infty$ *and* $\mathbb{E}\left[Y^4 \middle| Z = z, X = x\right]$ *is bounded.*

*d.* $\mathbb{E}\left[Y^2 \middle| Z = z, X = x\right]$ *is continuously differentiable.*

*e.* $\mathbb{E}\left[Y \middle| Z = z, X = x\right]$ *and* $f_{Z,X}(z,x)$ *are* $s+1$ *times continuously, boundedly differentiable.*

Under these conditions, we can show the following asymptotic optimality and rates.

THEOREM 6. *Under Assumptions 2, 3, 4, and 5, we have that*

$$\sqrt{nh_n}(R(z) - \overline{R}_n(z)) \xrightarrow{d} \mathcal{N}(0,\, \eta_z \kappa) \quad \forall z \in \mathcal{Z},$$
$$\sqrt{nh_n^3}(z^\star - \overline{z}_n) \xrightarrow{d} \mathcal{N}\left(0,\, \frac{\eta_{z^\star}\kappa'}{R''(z^\star)^2}\right)$$
$$\left(nh_n^3\right)(R(z^\star) - R(\overline{z}_n)) \xrightarrow{d} \frac{-\eta_{z^\star}\kappa'}{2R''(z^\star)}\chi_1^2,$$

*and, if also* $nh_n^{2s+1} \to 0$, *then*

$$\sqrt{nh_n}(R(z^\star) - \overline{R}_n(\overline{z}_n)) \xrightarrow{d} \mathcal{N}(0,\, \eta_{z^\star}\kappa),$$

*where* $\mathcal{N}(0,\, \sigma^2)$ *denotes a centered normal distribution with variance* $\sigma^2$, $\chi_1^2$ *denotes a chi-squared distribution with one degree of freedom, and* $\eta_z$, $\kappa$, $\kappa'$ *are constants defined as follows*

$$\eta_z = r(z)^2\ \mathbb{E}\left[\frac{\mathrm{Var}\left(Y \middle| Z = z, X\right)}{f_{Z|X}(z|X)}\right], \quad \kappa = \int \tilde{K}(z)^2 dp, \quad \kappa' = \int \tilde{K}'(z)^2 dz,$$

*where* $\tilde{K}(z) = \int K(z,x)dx$ *and* $f_{Z|X}(z|x) = f_{Z,X}(z,x)/f_X(x)$ *is the conditional density of* $Z$.

The proof is given in the appendix.

The main implication of Theorem 6 is that, under regularity conditions but without model specification, the non-parametric solution $\overline{z}_n$ has objective value that converges to optimal with rate of convergence $1/n$. Note that Assumption 3 implies that $s \geq k$ when $k \geq 3$ and $s \geq k+1$ when $k \leq 2$. This means that a kernel of order strictly greater than two, also known as a "bias-reducing" kernel (Hansen 2009), is necessary when $k \geq 2$.

### 4.3. A Test Statistic and Large Sample Theory

The impediment to verifying our hypothesis $H_0$ is that $R(z)$, $z^\star$, and $\hat{z}$ are all unknown; were they known, we would compute $\rho = R(z^\star) - R(\hat{z})$ and compare it to 0. Therefore, we must come up with an observable test statistic as a proxy to $\rho$. We do this by replacing the unknowns by our consistent estimates for them. We replace $R(z)$ and $z^\star$ by our non-parametric estimates $\overline{R}(z)$ as in (5) and $\overline{z}$ as in (6) and we replace $\hat{z}$ by the "estimate" $\hat{z}_n$. The resulting test statistic is

$$\rho_n = \overline{R}_n(\overline{z}_n) - \overline{R}_n(\hat{z}_n). \tag{7}$$

If $\rho_n$ is small, we have reason to believe that $\rho = 0$, whereas if $\rho_n$ is large, we would believe that $\rho > 0$. The question is where to draw the line.

THEOREM 7. *Suppose Assumptions 2, 3, 4, 5, and 1 hold. Let* $\Gamma = \frac{-\eta_{z^\star}\kappa'}{2R''(z^\star)}$. *Then,*

*i.* *under* $H_0$, $(nh_n^3)\,\rho_n \xrightarrow{d} \Gamma\chi_1^2$, *and*

*ii.* *under* $H_1$, $(nh_n^3)\,\rho_n \xrightarrow{d} \infty$.

The proof is given in the appendix.

Theorem 7 says that if we only reject $H_0$ when $\rho_n > n^{-1}h_n^{-3}\Gamma F_{\chi_1^2}^{-1}(1-\alpha)$ (where $F_{\chi_1^2}^{-1}$ is the chi-squared quantile function), then when $H_0$ is true we would only falsely reject $H_0$ at most $\alpha$ fraction of the time (asymptotically). On the other hand, if $H_0$ is false, then we would eventually reject it using such a procedure (a property known as *consistency* of a hypothesis test). The problem is that $\Gamma$ is unknown meaning that this exact procedure cannot be implemented in practice.

### 4.4. A Hypothesis Test

One way to implement a hypothesis is to estimate $\Gamma$ and replace the estimate into the results of Theorem 7. In particular, given any estimate $\hat{\Gamma}_n$ that converges in probability to $\Gamma$, we would have as an immediate consequence of Theorem 7 that $(nh_n^3)\,\hat{\Gamma}_n^{-1}\rho_n$ converges in distribution to $\chi_1^2$ under $H_0$ and to $\infty$ under $H_1$. This would give an implementable test. Non-parametric estimators for $\Gamma$, however, would tend to be convoluted and unwieldy, involving partial means of estimators of conditional variance and density as well as fragile estimates of second derivatives of partial means.

Instead, we use the bootstrap (Efron and Tibshirani 1993) and the following observation.

THEOREM 8. *Suppose Assumptions 2, 3, 4, and 5 hold. Let* $A_n = \overline{R}_n(\overline{p}_n) - \overline{R}_n(p^\star)$. *Then,* $(nh_n^3)\, A_n \xrightarrow{d} \Gamma\chi_1^2$. *Consequently,* $(nh_n^3)\,\mathbb{E}\,[A_n] \to \Gamma$.

The proof is given in the appendix.

So, to estimate $\Gamma$, we use a scaled estimate of the mean of $A_n$. In the spirit of Besbes et al. (2010), we use the bootstrap to achieve this. The fact that $A_n$ is asymptotically pivotal suggests that a bootstrap procedure could be particularly powerful (Horowitz 2001). This bootstrap procedure is also more attractive than convoluted kernel estimates of $\Gamma$ because it is less dependent on parameters and it deals more directly with the finite-sample distribution of $\rho_n$.

Given data $\mathcal{S}_n = \{(X_1, Y_1, Z_1), \dots, (X_1, Y_1, Z_1)\}$, $\hat{z}_n$, and a significance $\alpha \in (0,1)$, the full hypothesis test for $H_0$ proceeds as follows:

1. Compute $\overline{R}_n$ and $\overline{z}_n$ as in equations (5)-(6) based on data $\mathcal{S}_n$.
2. Compute $\rho_n$ as in equation (7).
3. Fix $B$ large. For $b = 1, \dots, B$:
    a. Draw $n$ samples with replacement from $\mathcal{S}_n$ to form the resampled dataset $\mathcal{S}_n^{(b)}$.
    b. Compute $\overline{R}_n^{(b)}$ and $\overline{p}_n^{(b)}$ as in equations (5)-(6) based on the data $\mathcal{S}_n^{(b)}$.
    c. Set $A_n^{(b)} = \overline{R}_n^{(b)}(\overline{z}_n^{(b)}) - \overline{R}_n^{(b)}(\overline{z}_n)$.
4. Let $\hat{\Gamma}_n = \frac{nh_n^3}{B} \sum_{b=B}^{n} A_n^{(b)}$.
5. Return p-value $p = 1 - F_{\chi_1^2}(\hat{\Gamma}_n^{-1} n h_n^3 \rho_n)$ and reject $H_0$ if $p < \alpha$.

To summarize, in this the section, we developed a hypothesis test for the hypothesis that an observational-data-driven decision $\hat{z}_n$ (or, its probability limit point), has objective value that is indistinguishable from optimal in Problem (1) to a statistically significant degree. We developed the test in the setting where Assumption 2 holds and used a non-parametric a three-step optimized partial-means kernel estimator to construct a test statistic for the hypothesis. In Section 6, we apply this test to both synthetic and real data, but first we review a parametric solution to which we will also compare.

## 5. A Parametric Solution

In the preceding section we developed a non-parametric decision-making algorithm that converged to optimal without requiring any model to be specified. Non-parametric approaches, however, can sometimes be unwieldy because their shapelessness makes them uninterpretable and they may be slow to converge. In fact, there is a growing body of work (Besbes et al. 2010, Besbes and Zeevi 2015) arguing that parametric models are often sufficient for managerial decision-making problems, as the model may need only fit well near the optimum or even just induce an acceptable decision, whether or not the model is correct. In particular, what matters is not model fit but objective performance. In this section, we present a parametric way to estimate $y(z)$ from observational data under Assumption 2 using a generalization of the propensity score. The intention is to be able to study the power of parametric approaches in the observational data setting and see if they suffice there as well for decision-making purposes.

The propensity score is a common matching metric used in the comparison of binary treatments in observational data (Rosenbaum and Rubin 1983). The (conventional) propensity score of a study subject is equal to the conditional probability of receiving a treatment (rather than control) given the subject's covariates $X$. If treatments are continuous, the generalized propensity score of a unit is defined as the conditional density of the unit receiving whatever treatment it did receive given the subject's covariates (Robins et al. 2000, Hirano and Imbens 2004, Imai and Van Dyk 2004).

Following Hirano and Imbens (2004), we define the generalized propensity score as $Q = f_{Z|X}(Z, X)$, assuming the conditional density $f_{Z|X}(z|x)$ exists. That is, we take the conditional density $f_{Z|X}(z|x)$, which is non-random but unknown, and plug in as values the random variables $Z$ and $X$. The key property of the generalized propensity score is that it is sufficient as a control for identifying the mean-response curve $y(z)$ as summarized below in an adaptation of a common result.

THEOREM 9. *Suppose Assumption 2 holds and that $f_{Z|X}(z|x)$ exists. Then*

$$y(z) = \mathbb{E}\left[y(z, f_{Z|X}(z, X))\right], \text{ where } y(z, q) = \mathbb{E}\left[Y \middle| Z = z, Q = q\right].$$

The proof is given in the appendix for completeness. The implication of Theorem 9 is that it is sufficient to control just for the univariate generalized propensity score rather than all of $X$.

Following Hirano and Imbens (2004), Theorem 9 motivates the following general strategy:

1. Regress $Z$ on $X$ by fitting a generalized linear model (GLM) in order to estimate $f_{Z|X}(z,x)$. I.e., choose $\hat{\beta}_n$, $\hat{\tau}_n$ by maximum likelihood estimation, given the parametric model
$$f_{Z|X}(z|x;\beta,\tau) = h(z,\tau)\exp\left(\frac{b(\beta_0+\beta^T x)T(z) - A(\beta_0+\beta^T x)}{d(\tau)}\right).$$
See McCullagh et al. (1989) for choices of $b$, $T$, $A$, $d$, $h$. For example, the choices $b(\mu)=\mu$, $T(z)=z$, $A(\mu)=\mu^2/2$, $d(\tau)=\tau^2$, and $h(z,\tau)=\frac{1}{\sqrt{2\pi\tau}}e^{-\frac{z^2}{2\tau^2}}$ lead to ordinary least squares (OLS). Other examples of GLMs include logistic regression, Poisson regression, Gamma regression, and loglinear regression.
2. Use the fitted GLM to impute generalized propensity scores, setting $\widehat{Q}_i = f_{Z|X}(Z_i|X_i;\hat{\beta}_n,\hat{\tau}_n)$.
3. Regress $Y$ on $Z$ and $\widehat{Q}$ based on the imputed data $\{(Z_i,Y_i,\widehat{Q}_i): i=1,\dots,n\}$ using another GLM (e.g., linear or logistic regression) to produce an estimate $\hat{y}_n(z,q)$ of $y(z,q)$. For example, we can fit $Y=b^{-1}(\alpha_0+\alpha_1 z+\alpha_2 q+\alpha_3 q^2+\epsilon)$ via link function $b$ (e.g., if a log-log demand model is appropriate as in many pricing problems, we regress $\log(Y)$ on $\log(Z)$ and $\widehat{Q}$).
4. Use these to estimate the mean-response curve and prescribe the decision that optimizes the estimated objective,
$$\hat{z}_n \in \arg\max_{z\in\mathcal{Z}}\left\{r(z)\times\tfrac{1}{n}\sum_{i=1}^n \hat{y}_n(z,\hat{f}_{Z|X}(z|X_i;\hat{\beta}_n,\hat{\tau}_n))\right\}.$$

The above procedure provides a flexible parametric framework for computing $\hat{z}_n$ from observational data under Assumption 2. When we apply it to examples with both real and synthetic data in Section 6, we find that it performs well and produces rewards that are often statistically indistinguishable from optimal.

## 6. Empirical Investigation

In this section, we use first use simulated data to investigate how our test performs in a controlled environment and then use our test to study the power and limits of both predictive and parametric

**Figure 3 Predictive demand and revenue in Example 3**

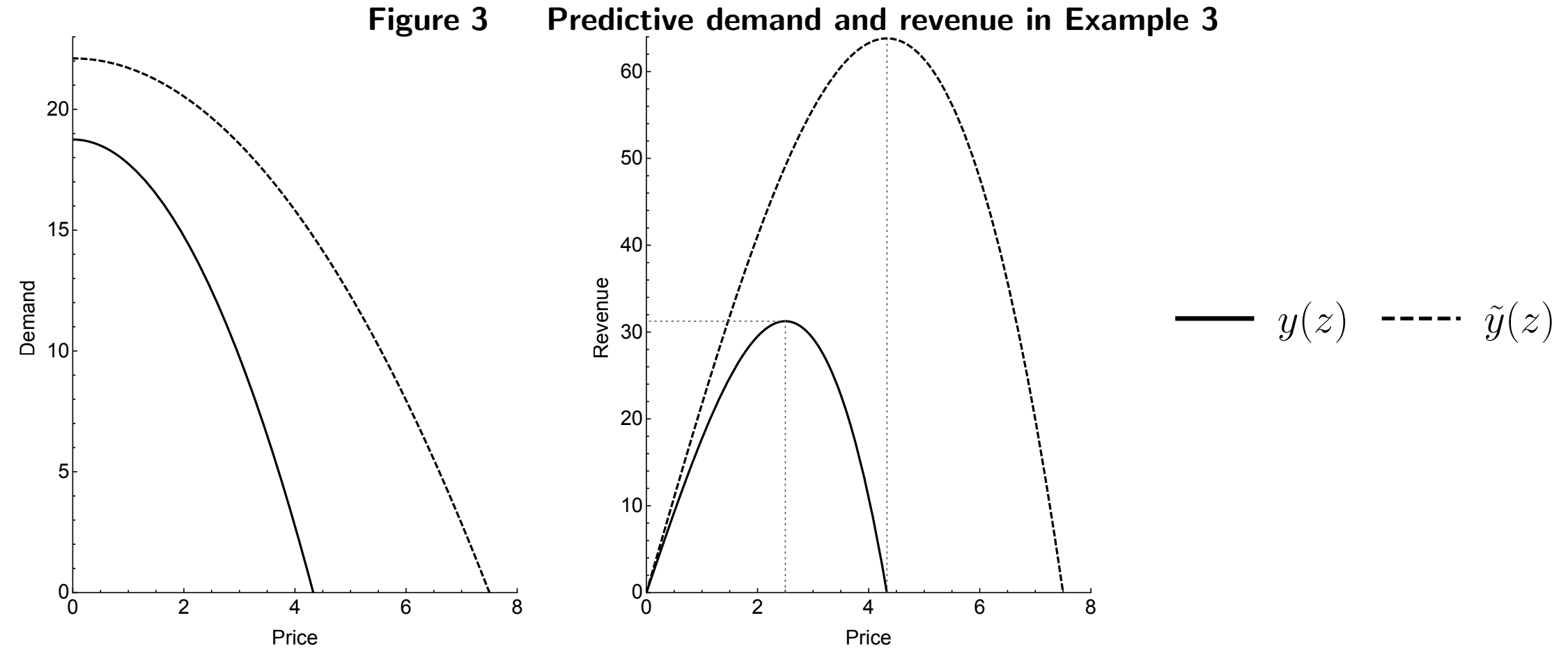

(a) Demand (b) Revenue

approaches to observational-data-driven decision-making in a real example. In the case of predictive approaches, we find that certain parametric ones sometimes yield good performance and sometimes not so good but are usually distinguishable from optimal, suggesting that predictive approaches are somewhat effective with both power and limits and, moreover, that the distinction between the true causal-effect objective $R(z)$ and the predictive surrogate $\tilde{R}(z)$ has real importance and practical relevance. In the case of parametric approaches, we find that they are almost always sufficient for good performance as long as it takes into account this distinction and correctly model the causal nature of $R(z)$ and the observational nature of the data, even if the parametric form may be misspecified.

Example 3 (Simulated Example). Consider a pricing instance of Problem (1) with procurement cost $c = 0$, potential prices $\mathcal{Z} = (0, \infty)$, and demand curve

$$Y(z) = 27.75 - z^2 + 6Xz - 9X^2 + V, \tag{8}$$

where $X \sim \mathcal{N}(0,1)$, $V \sim \mathcal{N}(0,\sigma^2)$ are normal noise, and prices historically set as $Z = 3X + W$ where $W \sim \mathcal{N}(0, \tau^2 = 15.1234)$.

The best predictor of $Y$, given an observation of $Z = z$, is

$$\mathbb{E}\left[Y \middle| Z = z\right] = 27.75 - z^2 + 6z\mathbb{E}\left[X \middle| 3X = z - W\right] - 9\mathbb{E}\left[X^2 \middle| 3X = z - W\right] + \mathbb{E}\left[V\right] = 22.108 - 0.393z^2,$$

**Figure 4 Comparing Predictive and Prescriptive Data-Driven Pricing Strategies**

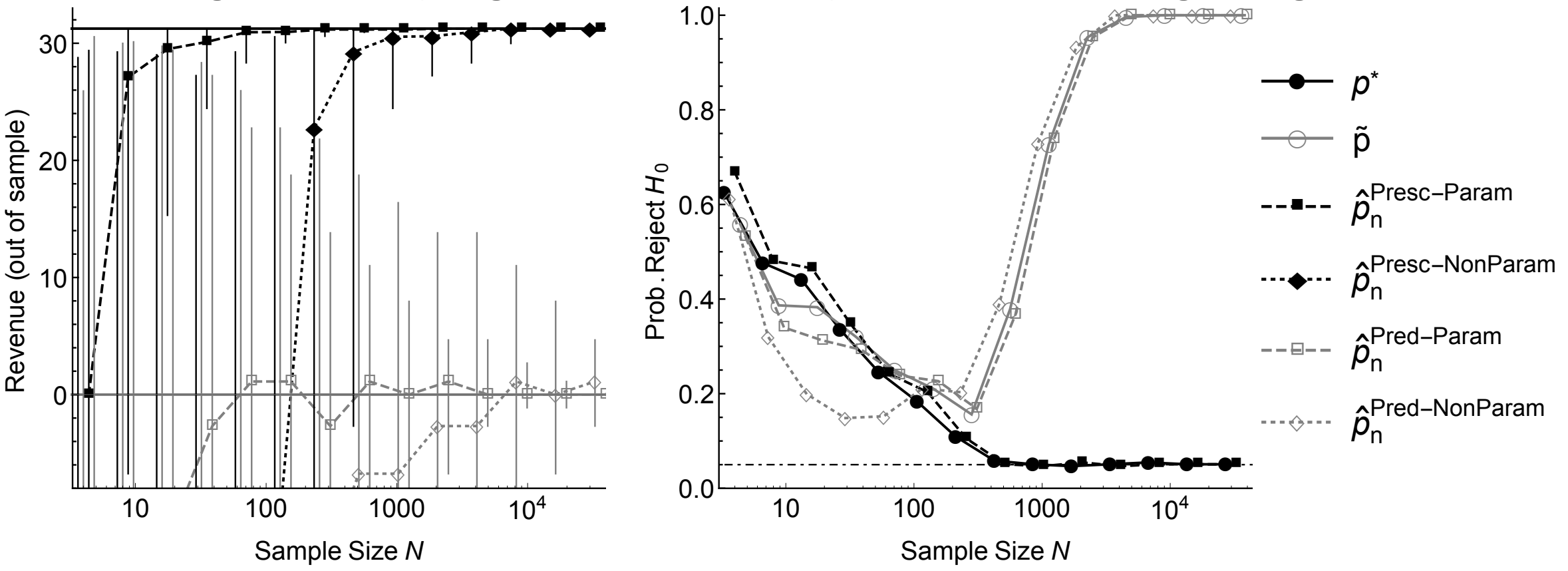

(a) Revenues (b) Test results

which we get by plugging in $P$ into (8) and recognizing $(X \mid 3X = z - W) \sim \mathcal{N}\left(\frac{3z}{9+\tau^2}, \frac{\tau^2}{9+\tau^2}\right)$. This is exactly the function we would arrive at if we used data to regress demand on price (e.g., by linear regression on $Z$, $Z^2$ or by non-parametric regression on $Z$). However, the expected demand when the price is set to $z$ is a different function,

$$\mathbb{E}\left[Y(z)\right] = 27.75 - z^2 + \mathbb{E}\left[X\right] - 9\mathbb{E}\left[X^2\right] + \mathbb{E}\left[V\right] = 18.75 - z^2,$$

which we get by taking the expectation of (8). We plot these two functions in Figure 3(a).

Now consider the price optimization problem. The true profit function, $R(z)$, is optimized at $z^\star = 2.5$ with a value of $R(z^\star) = 31.25$. On the other hand, a predictive approach optimizes $\tilde{R}(z)$, leading to the price $\tilde{z} = 4.330$, which leads to exactly $R(\tilde{z}) = 0$ profit under the true profit function. We plot $R$, $\tilde{R}$, $z^\star$, and $\tilde{z}$ in Figure 3(b).

We next compare four different observational-data-driven pricing: a prescriptive non-parametric approach $\hat{z}_n^{\text{Presc-NonParam}}$ based on $n$ observations of $(X, Y, Z)$, a prescriptive parametric approach $\hat{z}_n^{\text{Presc-Param}}$ based on $n$ observations of $(X, Y, Z)$, a predictive non-parametricapproach $\hat{z}_n^{\text{Pred-NonParam}}$ based on $n$ observations of $(Y, Z)$, and a predictive parametric approach $\hat{z}_n^{\text{Pred-Param}}$ based on $n$ observations of $(Y, Z)$. We also consider the (non-data-driven) true optimal price $z^\star$ of (1) and full-information predictive pricing strategy $\tilde{z}$ of (2). The prescriptive non-parametric strategy is as in (6) using a second order Gaussian kernel $K(u) = e^{-\frac{\|u\|_2^2}{2h_n^2}}$ and $h_n = 0.1 \times (n \log(n))^{-1/7}$, which satisfy

Assumption 3 with $s=2$, $k=1$. For the prescriptive parametric strategy, we follow our procedure from Section 5 using OLS linear regression of $P$ on $X$ for the GLM in step 1 and an OLS linear regression of $D$ on $P$ and $\log(\widehat{Q})$ in step 3. For the predictive non-parametric strategy, we use kernel regression to regress $r(Z)Y$ on $Z$ (using the same kernel and bandwidth $h_n = 2.5 \times (n\log(n))^{-1/7}$). Finally, for the predictive parametric strategy, we perform OLS linear regression of $Y$ on $Z$.

First, we consider the profit performance of each of these strategies. We plot the corresponding out-of-sample profits, $R(\hat{z}_n)$, along with optimal profit $R(z^\star)$, in Figure 4(a). The plot displays the median profit (center lines) and the $10^{\text{th}}$ and $90^{\text{th}}$ percentiles (vertical lines) over 256 replicate runs of each sample size. We see that the predictive approaches, by design, have very low revenues because of significant confounding. On the other hand, the prescriptive parametric approach offers significantly better out-of-sample performance than the non-parametric approach for small samples. In this example, the parametric approach is well-specified by design.

Next, we apply our hypothesis test for profit optimality. We plot the frequencies of rejecting a pricing strategy as significantly suboptimal at a significance of 0.05 in Figure 4(b). The plot displays the fraction of times the null hypothesis is rejected out of 256 replicate runs of each sample size. We see that with sufficient data, the test can distinguish those pricing strategies that generate suboptimal profits (i.e. solely predictive strategies) from those that cannot be distinguished from optimal for all prescriptive intents and purposes. In particular, it takes a few hundred data points before the test has the desired significance of 0.05 (i.e., $z^\star$ is rejected no more than 5% of the time).

Example 4 (Auto Loan Rate Optimization). Consider again the online auto loan dataset from Section 4.1.1. In Besbes et al. (2010), the authors consider whether a parametric model suffices for the problem of fixed pricing within various customer segments of loan applicants, defined in terms of three factors:

1. FICO score: (690, 715] (range 1) or (715, 740] (range 2),
2. Loan term in months: $\leq 36$ (class 1), (36, 48] (class 2), (48, 60] (class 3), or $> 60$ (class 4).
3. Season: first half of data (half 1) or second half (half 2).

**Table 4 Testing Revenue Optimality in the Auto Loan Rate Optimization Example**

| | | FICO range 1 (690, 715] | | FICO range 2 (715, 740] | | |
|---|---|---|---|---|---|---|
| | | Half 1 | Half 2 | Half 1 | Half 2 | |
| $n$ | | 1359 | 732 | 1386 | 781 | Term cl. 1 |
| $\rho_n$ ($p$-val) | Prescriptive, Param | 0.37 (0.15) | 0.21 (0.11) | 0.24 (0.49) | 0.23 (0.018*) | |
| | Predictive, Param | 0.86 (0.030*) | 0.50 (0.013*) | 0.25 (0.48) | 0.70 (< 0.001***) | |
| | Predictive, Non-Param | 1.91 (0.0012**) | 1.51 (< 0.001***) | 1.6 (0.07) | 1.35 (< 0.001***) | |
| $n$ | | 1394 | 832 | 1327 | 690 | Term cl. 2 |
| $\rho_n$ ($p$-val) | Prescriptive, Param | 0.23 (0.21) | 0.18 (0.073) | 0.28 (0.053) | 0.074 (0.67) | |
| | Predictive, Param | 0.87 (0.015*) | 0.24 (0.039*) | 0.35 (0.033*) | 0.051 (0.73) | |
| | Predictive, Non-Param | 1.6 (0.0011**) | 1.59 (< 0.001***) | 1.19 (< 0.001***) | 1.76 (0.040*) | |
| $n$ | | 4495 | 3147 | 3803 | 2865 | Term cl. 3 |
| $\rho_n$ ($p$-val) | Prescriptive, Param | 0.55 (0.32) | 0.26 (0.33) | 0.088 (0.061) | 1.4 (0.066) | |
| | Predictive, Param | 1.19 (0.14) | 0.22 (0.37) | 0.28 (< 0.001***) | 1.89 (0.034*) | |
| | Predictive, Non-Param | 1.19 (0.14) | 1.1 (0.046*) | 0.86 (< 0.001***) | 2.49 (0.015*) | |
| $n$ | | 2347 | 1506 | 1834 | 1206 | Term cl. 4 |
| $\rho_n$ ($p$-val) | Prescriptive, Param | 0.40 (0.0071**) | 0.0059 (0.63) | 1.86 (0.30) | 0.14 (0.46) | |
| | Predictive, Param | 0.27 (0.026*) | 0.045 (0.19) | 2.19 (0.26) | 0.31 (0.28) | |
| | Predictive, Non-Param | 1.5 (< 0.001***) | 1.54 (< 0.001***) | 2.92 (0.19) | 1.7 (0.012*) | |

* denotes reject $H_0$ at significance 0.05, ** at 0.01, and *** at 0.001. Gray denotes $p$-value $\geq 0.05$.

Customers with FICO scores outside of (690, 740] are not considered (see Besbes et al. 2010, for reasoning). Term classes 2 and 4 are not considered either, but we consider these here. The authors use a per-unit profit function $r(z) = z - 2\%$. Within each segment, the authors' approach is to estimate (either parametrically or non-parametrically) the conditional expectation of demand given price and to optimize per-unit profit times this conditional expectation. Using a test that compares the parametric and non-parametric approaches, they conclude that a parametric model suffices.

We consider the same problem again here, paying closer attention to the observational nature of the data. In Section 4.1.1 we argued that even within each segment, the data cannot be treated as experimental (i.e., satisfying Assumption 2 with respect to segment alone) and therefore that purely predictive approaches may not be estimating the true demand-price response function. We now use our hypothesis test to determine whether this distinction is moot from a profit-generated point of view. We also test whether our parametric prescriptive approach from Section 5 is successful. For the predictive approaches, we reproduce those in Besbes et al. (2010): kernel regression with the Gaussian kernel (non-parametric) and logistic regression (parametric). For our parametric prescriptive approach we fit a log-normal model for price via linear regression on $X$, i.e., $\big(\log(P)\big|X=x\big)\sim\mathcal{N}\left(\beta_0+\beta_1^T x,\,\sigma^2\right)$, and we fit a logistic regression for demand that is linear in price and quadratic in generalized propensity score, i.e., $\hat{d}(p,q)=\left(1+e^{-\alpha_0-\alpha_1 p-\alpha_2 q-\alpha_3 q^2}\right)^{-1}$.

We let $X$ consist of FICO score, the loan amount, the loan term, whether the car is new or used, whether the loan is refinancing, and if so what was the previous rate (otherwise 0). In our assessment, each of these covariates has direct impact on the interest rate quoted to applicants and each can arguably impact the decision of the applicant to accept any one rate. At the same time, this summarizes all relevant data provided and thus encapsulates all customer-specific information that could have gone into a rate quote decision. Therefore, we reason that Assumption 2 holds with respect to $X$, while it is likely to fail with respect to any subset of $X$.

We run the test within each of the 16 customer segments. In Table 4, we report the estimated suboptimality $\rho_n$ and its corresponding $p$-value according to our bootstrap procedure with $B=100$ draws.

The profits generated by the non-parametric predictive approach are rejected as suboptimal to a statistically significant degree at $p<0.05$ in 13 of 16 segments, and at $p<0.001$ in 7 segments. It is clear that this approach leaves much revenue on the table. The results for the parametric predictive approach are more mixed. The profits it generates are rejected as suboptimal to a statistically significant degree at $p<0.05$ in only 9 of 16 segments, and and at $p<0.001$ in only 2 segments. Moreover, the profit suboptimality estimate ($\rho_n$) is often small in magnitude. This shows that a good predictive approach can yield reasonable performance, but not without its limits. Indeed, a parametric approach that takes into account the observational nature of the data and has the true mean-response function (under Assumption 2) as its target (prescriptive, parametric) performs very well in this dataset. Our prescriptive parametric approach passes the test at $p\geq 0.05$ in all but 2 segments, in each of which, both predictive approaches also failed the test. Moreover, we see that the estimated suboptimality of our prescriptive parametric approach is smaller than that of the non-parametric predictive approach in all segments and than the parametric predictive approach in all but 3 segments. Averaging the estimated suboptimalities of each approach over all segments (weighted appropriately by the segment's $n$) and comparing, we see that on average, our prescriptive parametric approach recoups 70% (non-parametric) or 36% (parametric) of the total profits lost by using predictive approaches.

We note that in the same segments in which the analysis in Besbes et al. (2010) suggested that logistic regression on $P$ is sufficient for pricing purposes, our findings show that it is in fact insufficient. On the other hand, a parametric approach that addresses the observational nature of the data and the prescriptive nature of the problem seems to suffice in most cases, generating profits that the data cannot outright distinguish from optimal. In practice, it is known that parametric models, even if misspecified, can be helpful in extracting useful conclusions from smaller datasets. Our findings confirm this and, while refuting the evidence provided, agree with the final conclusion of Besbes et al. (2010), that parametric approaches work well for pricing.

## 7. Conclusions

We studied the data-driven optimization problem that arises from observational data. Noting that this problem is not generally well specified by the data, we considered another problem that is in fact the problem addressed by commonly used predictive approaches. While the two problems are different in their objective functions, we argued that predictive approaches will still work to the extent that their performance in the true problem is good. To quantify this, we focused on the pricing instance of the decision-making problem and proved strong performance guarantees for predictive approaches by leveraging the special structure of the pricing problem as well as the special characteristics of pricing data. In order to test in practice whether predictive approaches work, whether the distinction is practically a moot point, and whether a parametric approach that addresses the observational nature of the data works, we developed a hypothesis test for causal-effect objective optimality in the true optimization problem. The test was based on non-parametric estimates of this objective and its optimizer based on an ignorability assumption. Using consistency and asymptotic normality of these estimates, we established an asymptotic null distribution for a test statistic for our hypothesis. We used this null distribution and a bootstrap estimate of a nuisance parameter to develop the complete hypothesis test. Applying the test to the online auto loan dataset, we found that good predictive approaches can often yield good performance, but with limits, while approaches to observational-data-driven optimization often suffice in practice but only when they take into full account the data's observational nature and the problem's prescriptive rather than predictive nature.

## Endnotes

1. There is also an important stream of literature looking at pricing based on repeated experimentation (Bertsimas and Perakis 2006, Besbes and Zeevi 2009, Harrison et al. 2012), which does not consider pricing based on a *historical* dataset.

2. Note that our null hypothesis differs from the one considered by Besbes et al. (2010) in the definition of $R(z)$, that is, our $R(z) = \mathbb{E}\left[r(z)Y(z)\right]$ vs. Besbes et al. (2010)'s $\tilde{R}(z) = \mathbb{E}\left[r(z)Y \middle| Z = z\right]$.

## EC.1. Omitted Proofs

*Proof of Theorem 2* That $d(p)$ is linear and decreasing implies that $d(p) = d_0 - \lambda(p-c)$ with $\lambda > 0$. Hence, $R(p) = d_0(p-c) - \lambda(p-c)^2$, which is unimodal and uniquely maximized at $p^\star = c + d_0/(2\lambda)$ with value $R(p^\star) = d_0^2/(4\lambda)$. Let $\delta(p) = \mathbb{E}\left[\epsilon \middle| P = p\right]$, $\eta = d_0\gamma$. Then $|\delta(p)| \le \eta$. Note that

$$\mathbb{E}\left[D \middle| P = p\right] = \mathbb{E}\left[D(p) \middle| P = p\right] = \mathbb{E}\left[d(p) + \epsilon(p) \middle| P = p\right] = d(p) + \mathbb{E}\left[\epsilon \middle| P = p\right] = d_0 - \lambda(p-c) + \delta(p).$$

Hence, the theorem is trivial if $\eta = 0$ so let us assume $\eta > 0$.

Next we ask the question, what is the largest and smallest that the maximizer $\tilde{p}$ of $\tilde{R}(p)$ can be. By assumption, $|\delta(p)| \le \eta$ for all $p \in \mathcal{P}$. So, defining $\tilde{R}_{\delta_0}(p) := (p-c)\left(d_0 - \lambda(p-c) + \delta_0(p)\right)$, we are interested in

$$\tilde{p}_{\max} = \sup\left\{\sup\left(\arg\max_{p\in\mathcal{P}} \tilde{R}_{\delta_0}(p)\right) : |\delta_0(p)| \le \eta\right\}, \tag{EC.1}$$

$$\tilde{p}_{\min} = \inf\left\{\inf\left(\arg\max_{p\in\mathcal{P}} \tilde{R}_{\delta_0}(p)\right) : |\delta_0(p)| \le \eta\right\}, \tag{EC.2}$$

where we define $\sup(\varnothing) = -\infty$ and $\inf(\varnothing) = \infty$ without loss of generality because we assumed an optimizer $\tilde{p}$ exists for $\tilde{R}(p)$ so we are only interested in those functions $\delta(p)$ that induce a nonempty argmax. In what follows, define $\tilde{R}_+(p) = (p-c)\left(d_0 - \lambda(p-c) + \eta\right)$ and $\tilde{R}_-(p) = (p-c)\left(d_0 - \lambda(p-c) - \eta\right)$, which are both unimodal and uniquely maximized at $\tilde{p}_+ = c + (d_0+\eta)/(2\lambda)$ and $\tilde{p}_- = c + (d_0-\eta)/(2\lambda)$ respectively ($\tilde{p}_- < \tilde{p}_+$ because $\eta > 0$). Notice that $\tilde{R}_-(p) \le \tilde{R}_{\delta_0}(p) \le \tilde{R}_+(p)$ whenever $|\delta_0(p)| \le \eta$ with equality when $\delta_0(p) = \pm\eta$ is extremal.

First we argue that the bounds (EC.1)-(EC.2) are finite. For any $p \ge p' = c + \left(d_0 + \eta + 2\sqrt{\lambda + d_0\eta}\right)/(2\lambda)$ and $|\delta_0(p)| \le \eta$, since $\tilde{R}_+(p)$ is decreasing past $\tilde{p}_+$ and $p' \ge \tilde{p}_+$, we have that

$$\tilde{R}_{\delta_0}(p) \le \tilde{R}_+(p) \le \tilde{R}_+(p') = (d_0-\eta)^2/(4\lambda) - 1 < (d_0-\eta)^2/(4\lambda) = \tilde{R}_-(\tilde{p}_-) \le \tilde{R}_{\delta_0}(\tilde{p}_-).$$

Since $\tilde{p}_- \le p'$ we conclude that $\tilde{p}_{\max} \le p' < \infty$. Finally, since $\delta_0(p) = 0$ is feasible in (EC.1)-(EC.2), we have $c \le \tilde{p}_{\min} \le p^\star \le \tilde{p}_{\max} \le p'$.

Next we argue that in (EC.1) it is sufficient to consider functions $\delta_0(p)$ taking values in $\{-\eta,+\eta\}$ that are monotonic increasing, i.e. constant or step functions. Let $\delta_0(p)$ be feasible in (EC.1) and let $\tilde{p}_0=\sup\left\{\arg\max_{p\in\mathcal{P}}\tilde{R}_{\delta_0}(p)\right\}$. If at any $p_1\geq\tilde{p}_0$ we have $\delta_0(p_1)<\eta$, then increasing the value of $\delta_0(p_1)$ to $\eta$ can only increase the value of $\tilde{R}_{\delta_0}(p_1)$, which in turn may only increase the largest maximizer since $p_1\geq\tilde{p}_0$. Moreover, if at any $p_1<\tilde{p}_0$ we have $\delta_0(p_1)>-\eta$, then decreasing the value of $\delta_0(p_1)$ to $-\eta$ can only decrease the value of $\tilde{R}_{\delta_0}(p_1)$, which must already be at or below the maximal value and hence must leave the largest maximizer unchanged. The argument is unchanged even if $\tilde{p}_0$ is $\pm\infty$. A symmetric argument shows that in (EC.2) it is sufficient to consider functions $\delta_0(p)$ taking values in $\pm\eta$ that are monotonic decreasing.

Next we evaluate $\tilde{p}_{\max}$. Fix $\tilde{p}'=c+(\sqrt{d_0}+\sqrt{\eta})^2/(2\lambda)$ and let us consider the step function $\delta_{\max}(p)=\eta\mathbb{I}\left[p\geq\tilde{p}'\right]-\eta\mathbb{I}\left[p<\tilde{p}'\right]$. Since $\tilde{p}'>p_-$, $\tilde{R}_{\delta_{\max}}(p)$ is uniquely maximized on $(c,\tilde{p}')$ at $p_-$, with value $\tilde{R}_{\delta_{\max}}(p_-)=\tilde{R}_-(p_-)=(d_0-\eta)^2/(4\lambda)$. Since $\tilde{p}'>p_+$, $R_{\delta_{\max}}(p)$ is uniquely maximized on $[\tilde{p}',\infty)$ at $\tilde{p}'$, with value $\tilde{R}_{\delta_{\max}}(\tilde{p}')=\tilde{R}_+(\tilde{p}')=(d_0-\eta)^2/(4\lambda)$. Hence, $\arg\max_{p\in\mathcal{P}}\tilde{R}_{\delta_{\max}}(p)=\{p_-,\tilde{p}'\}$ and $\sup\{p_-,\tilde{p}'\}=\tilde{p}'$. Now we show that it is impossible to achieve a higher maximizer with $|\delta(p)|\leq\eta$, which would lead to $\tilde{p}_{\max}=\tilde{p}'$. By our previous argument we need only consider functions $\delta(p)$ taking values in $\pm\eta$ that are monotonic increasing. The constant functions taking values in $\pm\eta$ induce the maxima $\tilde{p}_-$ and $\tilde{p}_+$, both of which are smaller than $\tilde{p}'$. Next, consider any step function $\delta_0(p)=\eta\mathbb{I}\left[p\geq\tilde{p}_0\right]-\eta\mathbb{I}\left[p<\tilde{p}_0\right]$ with $\tilde{p}_0\neq p'$. If $\tilde{p}_0\leq\tilde{p}_+$ then, for any $p\neq p_+$, we have that $\tilde{R}_{\delta_0}(\tilde{p}_+)=\tilde{R}_+(\tilde{p}_+)>\tilde{R}_+(p)\geq\tilde{R}_{\delta_0}(p)$ since $\tilde{p}_+$ is the unique maximizer of $\tilde{R}_+(p)$; hence $\tilde{p}_+<p'$ is the unique maximum of $\tilde{R}_{\delta_0}(p)$. Consider $\tilde{p}_0>\tilde{p}_+$. Then, since $\tilde{p}_0>p_+>p_-$, $\tilde{R}_{\delta_{\max}}(p)$ is uniquely maximized on $(c,\tilde{p}_0)$ at $p_-$, with value $\tilde{R}_{\delta_0}(p_-)=\tilde{R}_-(p_-)=(d_0-\eta)^2/(4\lambda)$. Since $\tilde{p}_0>p_+$, $R_{\delta_{\max}}(p)$ is uniquely maximized on $[\tilde{p}_0,\infty)$ at $\tilde{p}_0$, with value $\tilde{R}_{\delta_0}(\tilde{p}_0)=\tilde{R}_+(\tilde{p}_0)$. If $\tilde{p}_0<\tilde{p}'$, then either of these potential maximizers are smaller than $p'$. If $\tilde{p}_0>\tilde{p}'$ then, since $\tilde{R}_+(p)$ is strictly decreasing past $p_+$ and $\tilde{p}'\geq p_+$, we have $\tilde{R}_{\delta_0}(\tilde{p}_0)=\tilde{R}_+(\tilde{p}_0)<\tilde{R}_+(\tilde{p}')=(d_0-\eta)^2/(4\lambda)=\tilde{R}_-(p_-)=\tilde{R}_{\delta_{\max}}(p_-)$. Hence $\tilde{p}_-<\tilde{p}'$ is the unique maximum of $\tilde{R}_{\delta_0}(p)$. When $\eta<d_0$, a symmetric argument applied to (EC.2) shows that $\tilde{p}_{\min}=c+\left(\sqrt{d_0}-\sqrt{\eta}\right)^2/(2\lambda)$. If $\eta\geq d_0$, the lower bound $\tilde{p}_{\min}=c$ is achieved by $\delta_-(p)$. Hence, $\tilde{p}_{\min}=c+\max\left\{0,\sqrt{d_0}-\sqrt{\eta}\right\}^2/(2\lambda)$.

To summarize, we conclude that since $\left|\mathbb{E}\left[\epsilon\middle|P\right]\right| \le \eta$, we must have

$$\tilde{p} \in [\tilde{p}_{\min}, \tilde{p}_{\max}] \quad \text{where} \quad \tilde{p}_{\min} = c + \frac{\max\left\{0, \sqrt{d_0} - \sqrt{\eta}\right\}^2}{2\lambda}, \tilde{p}_{\max} = c + \frac{(\sqrt{d_0} + \sqrt{\eta})^2}{2\lambda}.$$

Plugging these bounds into $R(p)$ we have

$$R(\tilde{p}_{\max}) = \frac{d_0^2 - 4d_0\eta - \eta^2 - 4\eta\sqrt{d_0\eta}}{4\lambda}, \qquad R(\tilde{p}_{\min}) = \begin{cases} \frac{d_0^2 - 4d_0\eta - \eta^2 + 4\eta\sqrt{d_0\eta}}{4\lambda} & \eta < d_0 \\ 0 & \eta \ge d_0 \end{cases}$$

Notice that if $\eta < d_0$ then $R(\tilde{p}_{\max}) = R(\tilde{p}_{\min}) - 2\eta\sqrt{d_0\eta}/\lambda \le R(\tilde{p}_{\min})$ and if $\eta \ge d_0$ then $R(\tilde{p}_{\max}) \le 0 = R(\tilde{p}_{\min})$. Therefore, $\min\{R(\tilde{p}_{\max}), R(\tilde{p}_{\min})\} = R(\tilde{p}_{\max})$. Since $R(p)$ is unimodal and $\tilde{p} \in [\tilde{p}_{\min}, \tilde{p}_{\max}]$, we have

$$R(\tilde{p}) \ge \min\{R(\tilde{p}_{\max}), R(\tilde{p}_{\min})\} = R(\tilde{p}_{\max}) = \frac{d_0^2 - 4d_0\eta - \eta^2 - 4\eta\sqrt{d_0\eta}}{4\lambda}.$$

Finally, using $R(p^\star) = d_0^2/(4\lambda)$,

$$\frac{R(\tilde{p})}{R(p^\star)} \ge 1 - 4\left(\frac{\eta}{d_0}\right) - 4\left(\frac{\eta}{d_0}\right)^{3/2} - \left(\frac{\eta}{d_0}\right)^2,$$

which is a univariate polynomial in $\sqrt{\eta/d_0}$. □

*Proof of Theorem 3* Recall that $\mathbb{E}[\epsilon] = 0$. That $\epsilon$ and $P$ are jointly normal implies that

$$\mathbb{E}\left[\epsilon\middle|P\right] = \zeta(P - \mu), \quad \text{where } \mu = \mathbb{E}[P],\ \zeta = \frac{\text{Cov}(\epsilon, P)}{\text{Var}(P)}.$$

By assumption of non-positive correlation, $\zeta \le 0$.

That $d(p)$ is linear and decreasing implies that $d(p) = d_0 - \lambda(p - c)$ with $\lambda > 0$. Hence, $R(p) = d_0(p - c) - \lambda(p - c)^2$, which is unimodal and uniquely maximized at $p^\star = c + d_0/(2\lambda)$ with value $R(p^\star) = d_0^2/(4\lambda)$. Also, recall from the proof of Theorem 2 that $\mathbb{E}\left[D\middle|P = p\right] = \mathbb{E}\left[D(p)\middle|P = p\right] = \mathbb{E}\left[d(p) + \epsilon(p)\middle|P = p\right] = d(p) + \mathbb{E}\left[\epsilon\middle|P = p\right]$ and hence $\tilde{R}(p) = d_0(p - c) - \lambda(p - c)^2 + \zeta(p - c)(p - \mu)$, which is unimodal and uniquely maximized at its critical point $\tilde{p} = (2\lambda c + d_0 - \zeta(c + \mu))/(2\lambda - 2\zeta)$ because it is feasible since $\lambda > 0 \ge \zeta$ and $(2\lambda - 2\zeta)(\tilde{p} - c) = d_0 - \zeta(\mu - c) \ge d_0 \ge 0$.

Plugging $\tilde{p}$ into $R(p)$ we get

$$R(\tilde{p}) = \frac{(d_0 - \zeta(\mu - c))(d_0(\lambda - 2\zeta) + \lambda\zeta(\mu - c))}{4(\lambda - \zeta)^2}.$$

Rearranging and using $R(p^\star)=d_0^2/(4\lambda)$, we have

$$\frac{R(\tilde p)}{R(p^\star)}=1-\left(\frac{\zeta}{\lambda-\zeta}\right)^2\left(\frac{d_0+\lambda(c-\mu)}{d_0}\right)^2=1-\left(\frac{\zeta}{\lambda-\zeta}\right)^2\left(\frac{\mathbb{E}\left[D\right]}{d_0}\right)^2,$$

where we plugged in $\mathbb{E}\left[D\right]=\mathbb{E}\left[D(p)\right]=\mathbb{E}\left[d_0-\lambda(P-c)\right]=d_0-\lambda(\mu-c)$. Moreover,

$$\sup_{\zeta\le 0}\left(\frac{\zeta}{\lambda-\zeta}\right)^2=\lim_{\zeta\to-\infty}\left(\frac{\zeta}{\lambda-\zeta}\right)^2=1.$$

Hence, we have the result in the statement of the theorem. □

*Proof of Theorem 4* We repeat the proof of Theorem 2 but note that if $E(z)$ is non-increasing then in (EC.1) it is sufficient to consider *constant* functions $\delta_0(p)$ taking values in $\{-\eta,+\eta\}$, which implies $\tilde p_m ax=\tilde p_+$. Therefore,

$$R(\tilde p_{\max})=\frac{d_0^2-\eta^2}{4\lambda},$$

since $\eta\le d_0$ is assumed. Since $R(p)$ is unimodal and $\tilde p\in[\tilde p_{\min},\,\tilde p_{\max}]$, we have

$$\begin{aligned}R(\tilde p)\ge\min\left\{R(\tilde p_{\max}),\,R(\tilde p_{\min})\right\}&=\min\left\{\frac{d_0^2-\eta^2}{4\lambda},\frac{d_0^2-4d_0\eta-\eta^2+4\eta\sqrt{d_0\eta}}{4\lambda}\right\}\\&\ge\frac{d_0^2-4d_0\eta-\eta^2+4\eta\sqrt{d_0\eta}}{4\lambda},\end{aligned}$$

since $\eta\le d_0$ is assumed. Finally, using $R(p^\star)=d_0^2/(4\lambda)$,

$$\frac{R(\tilde p)}{R(p^\star)}\ge 1-4\left(\frac{\eta}{d_0}\right)+4\left(\frac{\eta}{d_0}\right)^{3/2}-\left(\frac{\eta}{d_0}\right)^2,$$

completing the proof. □

*Proof of Theorem 6* Assumption 2 gives $R(p)=\mathbb{E}\left[\mathbb{E}\left[r(P)D\middle|P=p,X\right]\right]$, i.e. profit is given by taking a partial mean with $P=p$ fixed of the regression of $r(P)D$ on $P$ and $X$.

By Assumption 5 part i, there exists $\delta>0$ such that $[p^\star-\delta,\,p^\star+\delta]$ is contained inside the support of $P$. By Assumption 5 part i, $f_{P,X}(p,x)$ is bounded away from 0 on its support, and, by Assumption 5 part v, $\frac{\partial f_{P,X}(p,x)}{\partial x}$ is bounded. Hence,

$$\left|\frac{\partial}{\partial x}\log(f_{P,X}(p,x))\right|=\frac{\left|\frac{\partial f_{P,X}(p,x)}{\partial x}\right|}{f_{P,X}(p,x)}\le L<\infty$$

on the support of $(P,X)$. Therefore, we have that, for any $x$ and $|\psi|\le\delta$,

$$\log(f_{P,X}(p^\star+\psi,\,x))\le\log(f_{P,X}(p^\star,\,x))+L\delta,$$

and consequently,

$$\int\sup_{|\psi|\le\delta}f_{P,X}(p^\star+\psi,x)dx\le\int e^{L\delta}f_{P,X}(p^\star,\,x)dx<\infty.$$

By Assumption 5 part iii, there exists $M$ such that $\mathbb{E}\left[D^4\middle|P=p,X=x\right]\le M$ for all $p$, $x$. Combined, this yields

$$\int\sup_{|\psi|\le\delta}\left(1+\mathbb{E}\left[(r(P)D)^4\middle|P=p^\star+\psi,X=x\right]\right)f_{P,X}(p^\star+\psi,x)dx$$
$$\le(1+r(p^\star+\delta)^4M)\int\sup_{|\psi|\le\delta}f_{P,X}(p^\star+\psi,x)dx<\infty. \quad \text{(EC.3)}$$

To study our profit function estimator (5) for each fixed $p$, we employ Theorem 4.1 of Newey (1994) (henceforth, N in this proof), which provides convergence results for two-step kernel $m$-estimators, a specific case of which is our profit function estimator. We let the "trimming function" of N be $\tau(x)=\mathbb{I}\left[f_X(x)>0\right]$, an indicator for the compact support of $X$ (where its density is assumed bounded away from zero). Since our estimator (5) has $K$ both in the numerator and denominator, it is unchanged if we rescale $K$ by a positive constant. Similarly, the conditions of Assumption 3 remain unchanged. Hence, by Assumption 3 part i, without loss of generality we may assume $\int_{\mathbb{R}^{1+k}}K=1$. Then, Assumption K of N is satisfied with $\Delta=2$ by Assumption 3 parts i-iv. Assumption H of N is satisfied with $d=s+1$ by Assumption 5 part v. These constitute condition (ii) of N's Theorem 4.1. By Assumption 5 part i, there exists $c<p_{\max}<\infty$ such that $P\le p_{\max}$ almost surely. Let $r_{\max}=r(p_{\max})<\infty$. Then, by Assumption 5 part iii, $\mathbb{E}\left[(r(P)D)^4\right]\le r_{\max}\mathbb{E}\left[D^4\right]<\infty$ and $\mathbb{E}\left[(r(P)D)^4\middle|P=p,X=x\right]=r(p)\mathbb{E}\left[D^4\middle|P=p,X=x\right]$ are bounded. Combined with Assumption 5 part ii, we satisfy condition (i) of N's Theorem 4.1. Condition (iii) of N's Theorem 4.1 is satisfied by our choice of $\tau(\cdot)$ and by Assumption 5 part i. The first clause of condition (iv) of N's Theorem 4.1 is satisfied by our choice of $\tau(\cdot)$ and by Assumption 5 part ii. The second clause is satisfied by Assumption 5 parts iv-v combined with the fact that for any $m>0$, $\mathbb{E}\left[(r(P)D)^m\middle|P=p,X=x\right]=$

$r(p)^m\mathbb{E}\left[D^m \middle| P=p, X=x\right]$ and $r(p)^m$ is continuous. The third clause is satisfied by (EC.3). Since $X\in\mathbb{R}^k$ and $P\in\mathbb{R}$, condition (v) of N's Theorem 4.1 is satisfied by Assumption 3 parts v-vi. Applying N's Theorem 4.1 for each fixed $p\in\mathcal{P}$, we get

$$\sqrt{nh_n}(R(p)-\overline{R}_n(p)) \xrightarrow{d} \mathcal{N}(0,\,\eta_p\kappa)\ \ \forall p\in\mathcal{P},$$

where $\eta_p\kappa$ is a simplification of the asymptotic variance in eq. (14) in N.

To study the optimizer of our profit function estimator, we employ Flores (2005) (henceforth, F in this proof). Conditions (ii-vi) of F's Theorem 3 are satisfied in a similar way to the case of N's Theorem 4.1. Condition (i) of F's Theorem 3 is satisfied by Assumption 4 parts i and iii, condition (vii) by Assumption 5 part v, condition (viii) by Assumption 4 part iv, condition (ix) by Assumption 4 parts ii and iv, and finally condition (x) by Assumption 3 parts v-vi. Applying F's Theorem 3, we get

$$\sqrt{nh_n^3}(p^\star-\overline{p}_n) \xrightarrow{d} \mathcal{N}\left(0,\,\frac{\eta_{p^\star}\kappa'}{R''(p^\star)^2}\right), \tag{EC.4}$$

simplifying the asymptotic variance.

By Assumption 4 part iv and using Taylor's theorem to expand $R(p)$ around $p=p^\star$, there exists $p_n\in[\min(p^\star,\overline{p}_n),\,\max(p^\star,\overline{p}_n)]$ such that

$$R(\overline{p}_n)=R(p^\star)+R'(p^\star)(\overline{p}_n-p^\star)+\frac{1}{2}R''(p_n)(\overline{p}_n-p^\star)^2.$$

By first order optimality conditions, $R'(p^\star)=0$. Hence, rearranging, we have

$$R(p^\star)-R(\overline{p}_n)=-\frac{1}{2}R''(p_n)(\overline{p}_n-p^\star)^2. \tag{EC.5}$$

By continuous transformation of eq. (EC.4), we have

$$\left(nh_n^3\right)(\overline{p}_n-p^\star)^2 \xrightarrow{d} \frac{\eta_{p^\star}\kappa'}{R''(p^\star)^2}\chi_1^2. \tag{EC.6}$$

Eq. (EC.4) also implies $\overline{p}_n \xrightarrow{\mathbb{P}} p^\star$, which also implies $p_n \xrightarrow{\mathbb{P}} p^\star$ since $p_n$ is sandwiched between $\overline{p}_n$ and $p^\star$. Since $R''(p)$ is continuous, we also get by continuous transformation that

$$R''(p_n) \xrightarrow{\mathbb{P}} R''(p^\star). \tag{EC.7}$$

Combining eqs. (EC.5)-(EC.7), we get the desired result,

$$\left(nh_n^3\right)\left(R(p^\star)-R(\overline{p}_n)\right)\xrightarrow{d}\frac{-\eta_{p^\star}\kappa'}{2R''(p^\star)}\chi_1^2.$$

If $nh_n^{2s+1}\to 0$, then we also satisfy the conditions of F's Theorem 4 with equal bandwidths. Applying F's Theorem 4, we get

$$\sqrt{nh_n}(R(p^\star)-\overline{R}_n(\overline{p}_n))\xrightarrow{d}\mathcal{N}\left(0,\,\eta_{p^\star}\kappa\right),$$

simplifying the asymptotic variance. □

*Proof of Theorem 9* We begin by showing that $D(p)\perp\!\!\!\perp P\big|f_{P|X}(p,X)$. On the one hand we have

$$\begin{aligned}
f_{P|f_{P|X}(p,X)}(p|q) &= \mathbb{E}\left[\delta(P-p)\big|f_{P|X}(p,X)=q\right]\\
&= \mathbb{E}\left[\mathbb{E}\left[\delta(P-p)\big|f_{P|X}(p,X)=q,X\right]\big|f_{P|X}(p,X)=q\right]\\
&= \mathbb{E}\left[\mathbb{E}\left[\delta(P-p)\big|X\right]\big|f_{P|X}(p,X)=q\right]\\
&= \mathbb{E}\left[f_{P|X}(p,X)\big|f_{P|X}(p,X)=q\right]\\
&= q.
\end{aligned}$$

On the other hand, using Assumption 2, we have

$$\begin{aligned}
f_{P|f_{P|X}(p,X),D(p)}(p|q,d) &= \mathbb{E}\left[\delta(P-p)\big|f_{P|X}(p,X)=q,D(p)=d\right]\\
&= \mathbb{E}\left[\mathbb{E}\left[\delta(P-p)\big|f_{P|X}(p,X)=q,D(p)=d,X\right]\big|f_{P|X}(p,X)=q,D(p)=d\right]\\
&= \mathbb{E}\left[\mathbb{E}\left[\delta(P-p)\big|D(p)=d,X\right]\big|f_{P|X}(p,X)=q,D(p)=d\right]\\
&= \mathbb{E}\left[\mathbb{E}\left[\delta(P-p)\big|X\right]\big|f_{P|X}(p,X)=q,D(p)=d\right]\\
&= \mathbb{E}\left[f_{P|X}(p,X)\big|f_{P|X}(p,X)=q,D(p)=d\right]\\
&= q.
\end{aligned}$$

Equality between the two conditional probabilities implies the desired independence.

Using this independence and then plugging in $P=p$, we have

$$\mathbb{E}\left[D(p)\big|f_{P|X}(p,X)=q\right]=\mathbb{E}\left[D(p)\big|P=p,f_{P|X}(p,X)\right]=\mathbb{E}\left[D\big|P=p,Q=q\right]=d(p,q).$$

By iterated expectations, we get

$$\mathbb{E}[D(p)] = \mathbb{E}\left[\mathbb{E}\left[D(p)\middle|f_{P|X}(p,X)\right]\right] = \mathbb{E}\left[d(p, f_{P|X}(p,X))\right]$$

as desired. □

*Proof of Theorem 7* The proof borrows the outline of the proof of Theorem 2 of Besbes et al. (2010), but applied to our new testing case and causal estimators.

Decompose the test statistic $\rho_n$ into three terms:

$$\rho_n = \overline{R}_n(\overline{p}_n) - \overline{R}_n(\hat{p}_n) = A_n + B_n + C_n,$$

where

$$A_n = \overline{R}_n(\overline{p}_n) - \overline{R}_n(p^\star),$$

$$B_n = \overline{R}_n(p^\star) - \overline{R}_n(\hat{p}),$$

$$C_n = \overline{R}_n(\hat{p}) - \overline{R}_n(\hat{p}_n).$$

We begin by showing that $(nh_n^3)\,A_n \xrightarrow{d} \Gamma\chi_1^2$. By Assumption 3 part iii, we have that $\overline{R}_n(p)$ is twice continuously differentiable. Thus, using Taylor's theorem to expand $\overline{R}_n(p)$ around $p = \overline{p}_n$, we get that there exists $p_n \in [\min(p^\star, \overline{p}_n), \max(p^\star, \overline{p}_n)]$ such that

$$\overline{R}_n(p^\star) = \overline{R}_n(\overline{p}_n) + \overline{R}_n'(\overline{p}_n)(p^\star - \overline{p}_n) + \frac{1}{2}\overline{R}_n''(p_n)(p^\star - \overline{p}_n)^2.$$

By first order optimality conditions, $\overline{R}_n'(\overline{p}_n) = 0$. Hence, rearranging, we have

$$A_n = -\frac{1}{2}\overline{R}_n''(p_n)(p^\star - \overline{p}_n)^2. \tag{EC.8}$$

Next we show that $\overline{R}_n''(p_n) \xrightarrow{\mathbb{P}} R''(p^\star)$. Note that

$$\left|\overline{R}_n''(p_n) - R''(p^\star)\right| \leq \left|\overline{R}_n''(p_n) - R''(p_n)\right| + \left|R''(p_n) - R''(p^\star)\right|. \tag{EC.9}$$

As in the proof of Theorem 6, Assumptions 3, 4, and 5 imply the assumptions of Lemma 5.1 of Newey (1994) applied to $R''(p)$, which in turn yields the uniform convergence in probability of $\overline{R}_n''(p)$ over $\mathcal{P}$ since $\mathcal{P}$ is compact by Assumption 4 part i. Hence,

$$\left|\overline{R}_n''(p_n) - R''(p_n)\right| \leq \sup_{p\in\mathcal{P}}\left|\overline{R}_n''(p) - R''(p)\right| \xrightarrow{\mathbb{P}} 0. \tag{EC.10}$$

By Theorem 6, $\overline{p}_n \xrightarrow{\mathbb{P}} p^\star$. Because $p_n$ is sandwiched between $\overline{p}_n$ and $p^\star$, we also get $p_n \xrightarrow{\mathbb{P}} p^\star$. Since $R''(p)$ is continuous by Assumption 4 part iv, we have

$$|R''(p_n) - R''(p^\star)| \xrightarrow{\mathbb{P}} 0 \tag{EC.11}$$

by continuous transformation of the former. Combining eqs. (EC.9)-(EC.11), we get

$$\overline{R}_n''(p_n) \xrightarrow{\mathbb{P}} R''(p^\star). \tag{EC.12}$$

By continuous transformation of the result of Theorem 6 (eq. (EC.4)), we have

$$\left(nh_n^3\right)(\overline{p}_n - p^\star)^2 \xrightarrow{d} \frac{\eta_{p^\star}\kappa'}{R''(p^\star)^2}\chi_1^2. \tag{EC.13}$$

Combining eqs. (EC.8)-(EC.13), we get

$$\left(nh_n^3\right)A_n \xrightarrow{d} \frac{-\eta_{p^\star}\kappa'}{2R''(p^\star)}\chi_1^2 = \Gamma\chi_1^2. \tag{EC.14}$$

Next, we show that $(nh_n^3)\,C_n \xrightarrow{\mathbb{P}} 0$. By Assumption 3 part iii, we have that $\overline{R}_n(p)$ is twice continuously differentiable. Thus, using Taylor's theorem to expand $\overline{R}_n(p)$ around $p = \hat{p}$, we get that there exists $p_n' \in [\min(\hat{p}, \hat{p}_n), \max(\hat{p}, \hat{p}_n)]$ such that

$$\overline{R}_n(\hat{p}_n) = \overline{R}_n(\hat{p}) + \overline{R}_n'(\hat{p})(\hat{p}_n - \hat{p}) + \frac{1}{2}\overline{R}_n''(p_n')(\hat{p}_n - \hat{p})^2.$$

Rearranging, we have

$$\left(nh_n^3\right)C_n = -h_n^{3/2}\left(\sqrt{nh_n^3}\,\overline{R}_n'(\hat{p})\right)\left(\sqrt{n}(\hat{p}_n - \hat{p})\right) - \frac{1}{2}h_n^3\overline{R}_n''(p_n')\left(\sqrt{n}(\hat{p}_n - \hat{p})\right)^2. \tag{EC.15}$$

By Assumption 1, we have that

$$\sqrt{n}(\hat{p}_n - \hat{p}) = O_p(1), \text{ and hence also } \left(\sqrt{n}(\hat{p}_n - \hat{p})\right)^2 = O_p(1). \tag{EC.16}$$

Applying Theorem 4 of Newey (1994) we get the convergence in distribution of $\sqrt{nh_n^3}\left(\overline{R}_n'(p) - R'(p)\right)$ for any fixed $p$, including $\hat{p}$ and hence we have

$$\sqrt{nh_n^3}\,\overline{R}_n'(\hat{p}) = O_p(1). \tag{EC.17}$$

Next we show that $\overline{R}_n''(p_n') = O_p(1)$. Note that

$$\left|\overline{R}_n''(p_n') - R''(\hat{p})\right| \le \left|\overline{R}_n''(p_n') - R''(p_n')\right| + \left|R''(p_n') - R''(\hat{p})\right|. \tag{EC.18}$$

As before, $\overline{R}_n''(p)$ converges uniformly to $R''(p)$ in probability over $\mathcal{P}$ and so

$$\left|\overline{R}_n''(p_n') - R''(p_n')\right| \le \sup_{p\in\mathcal{P}} \left|\overline{R}_n''(p) - R''(p)\right| \xrightarrow{\mathbb{P}} 0. \tag{EC.19}$$

By Assumption 1, $\hat{p}_n \xrightarrow{\mathbb{P}} \hat{p}$. Because $p_n'$ is sandwiched between $\hat{p}_n$ and $\hat{p}$, we also get $p_n' \xrightarrow{\mathbb{P}} \hat{p}$. Since $R''(p)$ is continuous by Assumption 4 part iv, we have

$$\left|R''(p_n') - R''(\hat{p})\right| \xrightarrow{\mathbb{P}} 0 \tag{EC.20}$$

by continuous transformation of the former. Combining eqs. (EC.18)-(EC.20), we get

$$\overline{R}_n''(p_n') \xrightarrow{\mathbb{P}} R''(\hat{p}). \tag{EC.21}$$

Combining eqs. (EC.15)-(EC.21) gives $(nh_n^3)\,C_n = -h_n^{3/2}O_p(1) - h_n^3 O_p(1)$. Hence, because $h_n \to 0$, we get $(nh_n^3)\,C_n \xrightarrow{\mathbb{P}} 0$.

Finally, we treat $B_n$. Under $H_0$, $B_n = 0$ because Assumption 4 part ii (unique optimizer) and $H_0$ ($R(p^\star) = R(\hat{p})$) imply that $p^\star = \hat{p}$. Next, we show that under $H_1$, $(nh_n^3)\,B_n \xrightarrow{\mathbb{P}} \infty$. By applying the first results of Theorem 6 to each term, we have that $B_n \xrightarrow{\mathbb{P}} R(p^\star) - R(\hat{p})$. Since $k \ge 0$, Assumption 3 part vi implies $nh_n^5/\log(n) \to \infty$, which, since we also assume $h_n \to 0$, implies $nh_n^3 \to \infty$. Hence, since $R(p^\star) - R(\hat{R}) > 0$ under $H_1$, we have that $(nh_n^3)\,B_n \xrightarrow{\mathbb{P}} \infty$. □

*Proof of Theorem 8* Proven above. See eq. (EC.14). □